\documentclass[11pt]{amsart}
\usepackage[margin=1in]{geometry}
\usepackage[T1]{fontenc}
\usepackage[utf8]{inputenc} 
\usepackage{lmodern}
\usepackage{mathtools}   % includes amsmath
\usepackage{amssymb, amsfonts,mathrsfs}
\usepackage{amsthm}
\usepackage{bm}
\usepackage{graphicx}
\usepackage{booktabs, dcolumn}
\usepackage[tableposition=top]{caption}
\usepackage{tikz}
\usetikzlibrary{arrows.meta,calc,positioning}
\usepackage{xcolor}
\usepackage{hyperref}
\hypersetup{
    colorlinks=true,
    citecolor=blue,
    linkcolor=blue,
    filecolor=magenta,      
    urlcolor=cyan,
}
\usepackage{bookmark}   
\usepackage{enumerate}
\usepackage[shortlabels]{enumitem}
\numberwithin{equation}{section}
\newtheorem{theorem}{Theorem}[section]
\newtheorem{proposition}[theorem]{Proposition}
\newtheorem{lemma}[theorem]{Lemma}
\newtheorem{corollary}[theorem]{Corollary}
\newtheorem{claim}[theorem]{Claim}
\theoremstyle{definition}
\newtheorem{definition}[theorem]{Definition}
\theoremstyle{remark}
\newtheorem{remark}[theorem]{Remark}
\newcommand{\R}{\mathbb{R}}

\newcommand{\g}{\gamma}
\newcommand{\Z}{\mathbb{Z}}

\newcommand{\sph}{\mathbb{S}}

\newcommand{\Bry}{\operatorname{Bry}}
\newcommand{\Ric}{\operatorname{Ric}}
\newcommand{\Rm}{\operatorname{Rm}}
\begin{document}

\title{Geometry of symmetric 4d steady solitons}
\author{Daoyuan Han, Aprameya Girish Hebbar, Nata\v{s}a \v{S}e\v{s}um}
\date{}

\begin{abstract}
We prove that scalar curvature tends to zero at infinity on
complete, $\kappa$-noncollapsed, four-dimensional gradient steady Ricci
solitons with positive Ricci curvature, nonnegative sectional
curvature, and $\mathbb Z_2\times O(3)$ symmetry. 
\end{abstract}

%\keywords{}
%\subjclass{}
%\thanks{}

\maketitle
%\tableofcontents

\section{Introduction and main results}
\label{sec:introduction}

Gradient steady Ricci solitons are self-similar solutions to the Ricci
flow that evolve by diffeomorphisms. They arise as models for singularity
formation, and understanding their geometry is an important step toward
describing the possible singularities of the flow. A Riemannian manifold
$(M^n,g)$ together with a smooth function $f$ is called a gradient steady
Ricci soliton if
\begin{equation}\label{eq:steady-soliton-intro}
    \operatorname{Ric}_g+\nabla^2 f=0.
\end{equation}
If $\Phi_t$ is the family of diffeomorphisms generated by $\nabla f$,
then $g(t)=\Phi_t^*g$ is an eternal solution to the Ricci flow.

For dimensions $n\ge 3$, the Bryant soliton is a rotationally symmetric steady soliton \cite{Bry05}. It has positive curvature operator and linear scalar curvature decay. 
In dimension four, Lai's construction \cite{Lai20} provides a family of complete, $\kappa$-noncollapsed steady solitons with positive curvature operator and $\mathbb{Z}_2\times O(3)$ symmetry. 
Their asymptotic geometry differs from that of the Bryant soliton and this motivates the further study of $4d$ steady solitons. 
 
In \cite{HS26}, the second and the third named authors study complete
$4d$ steady solitons with positive Ricci curvature, nonnegative sectional curvature, a critical point of the potential, and $\kappa$-noncollapsing. In the non-Bryant case, they identified two distinguished edges, denoted by $\mathcal{G}=\mathcal{G}_{-}\cup\mathcal{G}_{+}$, and proved that the scalar curvature $R(x)$ satisfies 
\[
    R(x)\leq \frac{C}{d_g(x,\mathcal{G})},
    \qquad x\in M\setminus\mathcal{G}.
\]
A central question left open by this analysis was whether scalar curvature necessarily tends to zero at infinity.

In this paper, we prove that scalar curvature tends uniformly to zero
at infinity for $4d$ steady solitons with $\Z_2\times O(3)$ symmetry under the assumptions
stated below. Our argument combines the structural results of
\cite{HS26} with the mode analysis developed by Ma, Mahmoudian, and
the third named author in \cite{MMS23}. We carry out this analysis without
assuming scalar curvature decay and use reflection symmetry to exclude
the positive-mode alternative. 

From now on we normalize the metric so that $R+|\nabla f|^2=1.$ Our first main result of the paper is the following.
\begin{theorem}[Scalar curvature decay]
\label{thm:scalar-decay}
Let $(M^4,g,f)$ be a complete, connected gradient steady Ricci soliton
satisfying
$
    \operatorname{Ric}_g>0$, 
   $ \sec_g\geq 0,
$
and suppose that $(M,g)$ is $\kappa$-noncollapsed on all scales.
Assume that $f$ has a critical point $o$ and that $(M,g)$ is
$\mathbb{Z}_2\times O(3)$-symmetric in the sense of
\cite[Definition~2.1]{Lai20}. Then, 
\[
    \lim_{d_g(x,o)\to\infty}R(x)=0.
\]
If $(M,g)$ is not isometric to the $4d$ Bryant soliton, then the following conclusions hold.
\begin{enumerate}
    \item    The scalar curvature $R(x)$ satisfies 
    \[
        \limsup_{d_g(x,o)\to\infty}
        \frac{R(x)d_g(x,o)}{\log d_g(x,o)}
        =1.
    \]

    \item Let $\mathcal{G}=\mathcal{G}_{-}\cup\mathcal{G}_{+}$ be the union of the two edges constructed in \cite[Theorem~1.4]{HS26}. Then
    \[
        \limsup_{d_g(x,\mathcal{G})\to\infty}
        \frac{R(x)d_g(x,\mathcal{G})^2}
             {\log d_g(x,\mathcal{G})}
        =8.
    \]
    Consequently, for every $\sigma\in(0,2)$ there exists $C_\sigma>0$ such that
    \[
        R(x)\leq
        \frac{C_\sigma}
             {\bigl(1+d_g(x,\mathcal{G})\bigr)^\sigma},
        \qquad x\in M.
    \]
\end{enumerate}
\end{theorem}

In particular, the non-Bryant solitons in Theorem~\ref{thm:scalar-decay} do not satisfy a uniform linear scalar curvature decay estimate. The curvature-decay dichotomy conjecture of Munteanu--Sung--Wang \cite{MSW19}, recalled in \cite[Conjecture~1.1]{ChanZhu22}, is false in $4d$. 

Having established curvature decay, we study the geometry of
these $4d$ solitons. We set 
$\rho:=f(o)-f,
    \Sigma_s:=\{\rho=s\}$, so that $\nabla^2\rho=\operatorname{Ric}_g$, and $\rho(x)\to \infty$ as $d(x,o)\to \infty$.
\begin{theorem}[Asymptotic geometry]
\label{thm:asymptotic-geometry}
Let $(M^4,g,f)$ satisfy the hypotheses of
Theorem~\ref{thm:scalar-decay}, and assume that $(M,g)$ is not
isometric to the four-dimensional Bryant soliton. Then the following
hold.
\begin{enumerate}
    \item The asymptotic cone of $(M,g)$ is a ray. More precisely,
    \[
        (M,\lambda g,o)
        \xrightarrow[\lambda\downarrow 0]{\mathrm{pGH}}
        \bigl([0,\infty),d_{\mathbb{R}},0\bigr).
    \]

    \item The intrinsic diameters of $\Sigma_s=\{\rho=s\}$ satisfy
    \[
        \operatorname{diam}
        \bigl(\Sigma_s,g|_{\Sigma_s}\bigr)
        =4\bigl(1+o(1)\bigr)\sqrt{s\log s}
        \qquad\text{as }s\to\infty.
    \]
    Moreover,
    \[
        \operatorname{Vol}_g B_g(o,r)
        =\frac{128\pi}{15}\bigl(1+o(1)\bigr)
          r^{5/2}\sqrt{\log r}
        \qquad\text{as }r\to\infty.
    \]

    \item  Fix a reference point $q\in\rho^{-1}(s_0)$ whose integral curve
$q_s\in\Sigma_s$ of $\nabla\rho/|\nabla\rho|^2$, with
$q_{s_0}=q$, lies on the spherical orbit of maximal radius
in $\Sigma_s$ for every $s\ge s_0$. Write the induced metric on $\Sigma_s$ as
    $
        g|_{\Sigma_s}
        =dz^2+F(z,s)^2g_{\mathbb{S}^2},
    $
    where $g_{\mathbb{S}^2}$ is the round metric on $\sph^2$ with sectional curvature $1$ and $z$ is
    signed distance, measured from the spherical orbit through $q_s$. 
    For every fixed $L>0$, uniformly on $|z|\leq L\sqrt{s}$, \[
        F(z,s)
        =\sqrt{2s}
         -\frac{z^2-2s}{4\sqrt{2s}\log s}
         +o\left(\frac{\sqrt{s}}{\log s}\right)
        \qquad\text{as }s\to\infty.
    \]
\end{enumerate}
\end{theorem}

Next, we study the stability of these solitons under the Ricci flow.
Deruelle \cite{Deruelle} established a framework in which strict linear
stability implies dynamical stability for complete noncompact steady
solitons satisfying certain curvature assumptions. It is also proved that the Bryant soliton is dynamically stable. Theorems~\ref{thm:scalar-decay} and ~\ref{thm:asymptotic-geometry} allow us to apply the results in  \cite{Deruelle}. 
\begin{theorem}[Dynamical stability]
\label{thm:dynamical-stability}
Let $(M^4,g,f)$ satisfy the hypotheses of Theorem~\ref{thm:scalar-decay}. Then $(M,g)$ is dynamically stable.

More precisely, let $g_{\mathrm{in}}$ be a smooth metric on $M$ such
that
\[
    I:=
    \left(
        \int_M
        |g_{\mathrm{in}}-g|_g^2 e^{-f}\,d\mu_g
    \right)^{1/2}
    <\infty.
\]
There exists $\varepsilon=\varepsilon(g,I)>0$ such that, if
$
    \|g_{\mathrm{in}}-g\|_{L^\infty(g)}<\varepsilon,
$
then the modified Ricci--DeTurck flow with background
$(g,\nabla^g\rho)$ and initial metric $g_{\mathrm{in}}$ exists for
all $t\geq 0$ and converges exponentially to $g$ in every global
$C^k(g)$ norm. Consequently, the corresponding Ricci flow is immortal and converges
to $g$ modulo time-dependent diffeomorphisms.
\end{theorem}

The following corollary is immediate. 
\begin{corollary}
The above results apply to $\mathbb Z_2\times O(3)$-symmetric steady solitons constructed by Lai \cite{Lai20}. 
\end{corollary}

Recently, Fang \cite{Fang2026} proved that every complete, connected,
$4d$ gradient steady Ricci soliton with $\Rm>0$ that is $\kappa$-noncollapsed on all scales admits an effective, isometric $O(3)$ action with principal orbits diffeomorphic to $\mathbb{S}^2$, and the potential function has a critical point. In this setting, the
remaining hypothesis of Theorem~\ref{thm:scalar-decay} is an appropriate $\Z_2$-symmetry. More recently, Bamler and Lai \cite{BL26} develop a framework for deriving asymptotics for ancient asymptotically cylindrical Ricci flows. Our proof does not rely on their results and extends the methods of \cite{MMS23}, \cite{ABDS22}. 

\medskip
\noindent \textbf{Outline of the paper. }
Suppose $(M^4,g,f)$ is an $O(3)$-symmetric steady soliton with $\sec\geq 0,\Ric>0$, and $\nabla f(o)=0$. In \cite{MMS23}, the authors establish the dichotomy of modes under the additional assumption that $R\to 0$ at infinity, and prove that the dominance of the neutral mode gives fine asymptotics that resemble those of Perelman's ancient Type II solution.

In sections ~\ref{sec:dichotomy} and \ref{sec:neutral-mode}, after using the structural properties in \cite{HS26}, we closely follow the work of \cite{MMS23} and ensure that each step goes through without the additional assumption. In particular, we analyze the nonlinear parabolic PDE satisfied by $F(z,s)$ and show the dichotomy of modes (i.e., either the neutral mode or the positive modes dominate). Further, if the neutral mode dominates, we obtain $R\to 0$ at infinity as a consequence. 

In section~\ref{sec:positive-exclusion}, we assume $\Z_2$-symmetry and show that the positive mode can never dominate. The argument is by contradiction. Suppose that the positive modes dominate. We first use ideas inspired by \cite{ABDS22} to obtain estimates on the slope function $F_z(z,s)$. The proof uses a barrier which is constructed in Appendix~\ref{sec:appendix-A}. We then show that  $\operatorname{diam}(\Sigma_s)\geq cs$, hence obtain that $R\nrightarrow 0$ at infinity. Then, we use a logarithmic averaging argument to obtain a contradiction by analyzing how fast the slope decays along the level set $\Sigma_s$. 

Combining the work of sections ~\ref{sec:dichotomy}, \ref{sec:neutral-mode}, and ~\ref{sec:positive-exclusion}, under $\Z_2$-symmetry, it follows that $R\to 0$ at infinity and the neutral mode dominates. Using \cite{MMS23} and the asymptotics obtained there, we derive several properties of the soliton such as rate of curvature decay, diameter and volume growth. Finally combining these with the work of  \cite{Deruelle}, we obtain linear and dynamic stability for $(M,g)$. In Appendix~\ref{sec:appendix-B}, we recall a Sobolev inequality and an interpolation inequality used in the paper. 

\noindent \textbf{Statement on the use of AI. }
AI systems were used in finding the barrier functions $P_a$ in Section~\ref{sec:intermediate-region-slope} and Appendix~\ref{sec:appendix-A}, the barriers for the argument in Lemma~\ref{lem:Ric-lower-bound}, and the contradiction argument in Section~\ref{sec:coupling-contradiction}. 
The final mathematical arguments were simplified, written, and verified by the authors. 
The authors take full responsibility for the correctness and content of the paper. 

\noindent \textbf{Acknowledgements. }The authors would like to thank Zilu Ma for discussions.  N.S. gratefully acknowledges support from the National Science Foundation under grant DMS-2505574.

\section{Setting and Preliminaries}
\label{sec:setting}
We assume throughout the paper that $(M^4,g,f,o)$ is a complete gradient steady Ricci soliton satisfying the assumptions (A1)--(A3) of \cite{MMS23}: 
\begin{itemize}
\item [(A1)] $\Ric>0$, 
\item[(A2)] $\left(M^4, g\right)$ admits an ${O}(3)$-symmetry, that is, ${O}(3) \subset \operatorname{Isom}\left(M^4, g\right)$, 
\item[(A3)]  $f$ has a critical point $o$,
\end{itemize}
and $(M^4,g)$ has $\sec\geq 0$, is $\kappa$-noncollapsed in the sense of Perelman. We use the convention
\[
 \Ric_g+\nabla^2 f=0.
\]
This differs from \cite{MMS23} by a sign. The four-dimensional Bryant soliton has scalar curvature tending to zero at infinity and asymptotic cone equal to a ray. We therefore restrict attention to the non-Bryant case and assume that the tangent flow at infinity of $(M,g)$ is $\mathbb S^2\times\mathbb R^2$ (see the discussion in \cite[Section~1.1]{HS26}). Beginning in section~\ref{sec:positive-exclusion}, we will assume $\Z_2\times O(3)$ symmetry as defined in \cite[Definition~2.1]{Lai20}. 

Set $\rho:=f(o)-f\geq 0, \Sigma_s:=\{\rho=s\}.$ After scaling the metric, we have
\begin{equation}\label{eq:normalization}
 \nabla^2\rho=\Ric_g,\qquad R+|\nabla\rho|^2=1,
 \qquad \nabla R=-2\Ric_g(\nabla\rho,\cdot).
\end{equation}
Here and below, $R=R_g$, and $ n:=\frac{\nabla\rho}{|\nabla\rho|}$ denotes the outward unit normal to the level sets of $\rho$. Fix a large $s_0$ and put $\Sigma=\Sigma_{s_0}$. Define $\chi_s:\Sigma\to \Sigma_s$, $s\geq s_0$, 
\[
 \partial_s\chi_s=\frac{\nabla\rho}{|\nabla\rho|^2}\circ\chi_s,
 \qquad \chi_{s_0}=\mathrm{id}_\Sigma,\qquad
 \bar g_s=\chi_s^*(g|_{\Sigma_s}).
\]
Then $\chi_s(\Sigma)=\Sigma_s$. Let $\Phi_t,t\in \R$ be generated by $\nabla f=-\nabla\rho$, so that $g(t)=\Phi_t^*g,t\in \R$ solves Ricci flow. We refer to \cite{HS26} for further background and preliminary results. 
\begin{lemma}
\label{lem:exterior}
For each fixed $s_0>0$, there is $c_0>0$ such that
\[
 c_0\le |\nabla\rho|^2\le1\quad\hbox{on }\{\rho\ge s_0\}.
\]
Moreover, $\rho$ is proper, its positive level sets are compact and diffeomorphic to $\mathbb S^3$, and there exist positive constants $c_1,c_2$ such that for all $x\in M$, 
$$d(o,x)\geq \rho(x)\geq c_1d(o,x)-c_2.$$
\end{lemma}
\begin{proof}
The properness, comparison to $d_g(\cdot,o)$, and topology of level sets of $\rho$ follow from \cite[Lemma 2.1]{HS26} and the discussion after it. Since $s\mapsto R(\chi_s(q))$ is decreasing for each $q\in \Sigma$, we have $1\geq |\nabla \rho|^2=1-R\geq 1-\max_{\rho=s_0}R>0$. 
\end{proof}
We will be using the properties of the soliton proved in \cite[Section 1-5]{HS26}. On $M$, we have 
$$|\nabla^k\Rm|\le C_kR^{1+k/2},$$
see \cite{CMZ25} or \cite[Lemma~2.8(i)]{HS26}.  
\begin{definition}
\label{def:edges}
Let $\mathcal G_\pm:[0,\infty)\to M$ be the two edges in $M$, with $\mathcal G_\pm(0)=o$, $\mathcal{G}_{\pm}\cap \Sigma=\{x_\pm\}$ \cite[Definition~4.11]{HS26} and for $t>0$, 
\[
 \mathcal G_\pm'(t)=\frac{\nabla\rho}{|\nabla\rho|}(\mathcal G_\pm(t)).
\]
\end{definition}
We use $\nabla$ and $\Delta=\operatorname{tr}_g\nabla^2$ for the ambient connection and rough Laplacian.  The dimension-reduction and neck-cap results of \cite{HS26} give the structural results used in \cite[Sections 2 and 3.1]{MMS23},  with the additional remark that when the scalar curvature does not converge to zero, the level sets at the tips converge to a hypersurface in $\Bry^3\times \R$. 
\begin{proposition}[\cite{CMZ25,HS26}]
\label{prop:structure}
For sufficiently large $s$, $\Sigma_s$ has two tips $p_s^-$ and $p_s^+$. 
Moreover, the following hold. 
\begin{enumerate}
\item If $x_j\in M$ and $d(x_j,o)\to \infty$, then after passing to a subsequence $(M,R(x_j)g,x_j)$ converges smoothly to either $\Bry^3\times \R$ or $\sph^2\times \R^2$. 
    \item Let $s_j\to\infty$. Then after passing to a subsequence, $(M,R(p_{s_j}^\pm)g,p_{s_j}^\pm)$ converges to $\Bry^3\times \R$ in the smooth Cheeger--Gromov sense with the base point at the tip. 
    \item If $R(p_{s_j}^{\pm})\to 0$, then $(\Sigma_{s_j},R(p_{s_j}^\pm)g|_{\Sigma_{s_j}},p_{s_j}^\pm)$ converges to $(\Bry^3,g_{\Bry^3},x_{tip})$ in the smooth Cheeger--Gromov sense. 
    \item If instead $R(p_{s_j}^{\pm})\to\ell >0$, then $(\Sigma_{s_j},g|_{\Sigma_{s_j}},p_{s_j}^\pm)$ converges to $(S_h,g_{\Bry^3}+dz^2|_{S_h},(x_{tip},0))$ in the smooth Cheeger--Gromov sense, where 
\begin{equation}\label{eq:Sh}
 S_h=\{(y,t)\in\mathrm{Bry}^3\times\mathbb R:At+h(y)=0\},
 \qquad A\ne0,
\end{equation}
Here $h$ is a Bryant potential with $\Ric_{g_{\Bry^3}}+\nabla^2h=0$, $h(x_{tip})=0$, and $R_{g_{\Bry^3}}(x_{tip})=\ell, A^2=1-\ell.$ 

\end{enumerate}
\end{proposition}
\begin{proof}
All of these results are proved in \cite{CMZ25} and \cite[Lemma 2.8, Lemmas~3.1--3.2, Theorems~3.8 and~4.7]{HS26}. 
\end{proof}

\begin{lemma}
\label{lem:poles}
The tips $p_s^\pm$  satisfy $sR(p_s^{\pm})\to\infty$ as $s\to \infty.$ 
The two tips are the poles of the rotationally symmetric level sets. 
Moreover, $p_s^\pm =\chi_s x_\pm$ for all sufficiently large $s$. 
\end{lemma}
\begin{proof}
By \cite[Theorem~3.8]{HS26}, for all sufficiently large $s$, $\Sigma_s$ has exactly two tips $p_s^\pm$. 
The isometric action of $SO(3)$ preserves the set $\{p_s^+,p_s^-\}$. 
Since $SO(3)$ is connected and the map $SO(3)\to \{p_s^+,p_s^-\}$ given by $g\mapsto g\cdot p_s^+$ is continuous, it follows that the action preserves $p_s^+$. Similarly, the action also preserves $p_s^-$, hence fixes each tip. 
We conclude that the tips are the two poles of the rotationally symmetric level set.

By \cite[Lemma~2.6]{MMS23}, the function $\rho$ is invariant under the isometric action of $SO(3)$, hence $\chi_s$ is also invariant under that action. This proves that $p_s^\pm =\chi_s x_\pm$ for all sufficiently large $s$. 

We next prove that $sR(p_s^{\pm})\to\infty$ as $s\to \infty$. Fix either sign $\pm$ and put
$a(s)=R(p_s^\pm)$. 
If $a(s)\to L>0$ as $s\to \infty$, we are done. Suppose on the other hand that $a(s)\to 0$ as $s\to \infty$. 
Using the soliton identity $ \Delta R+\langle\nabla\rho,\nabla R\rangle
 =-2|\Ric|^2$, it follows that 
\[
 \frac{d}{ds}\frac1{a(s)}
 =
 \frac{\Delta R+2|\Ric|^2}
      {(1-a(s))a(s)^2}(p_s^\pm).
\]
Using Proposition~\ref{prop:structure}, it follows that $\frac{d}{ds}a(s)^{-1}\to 0.$
Integration yields $a(s)^{-1}=o(s)$, proving $sa(s)\to\infty$.
\end{proof}
The next proposition analyzes the behavior of intrinsic distance and the proof doesn't involve dimension reduction or symmetry assumptions. It holds for
any complete steady soliton with positive Ricci curvature and a critical point of the potential.
\begin{proposition}
\label{prop:zero-tip-sublinear}
Let $p, q\in \Sigma,p\neq q$ such that $\lim_{s\to \infty} R(\chi_sq)=0$. If $\lim_{s\to \infty} R(\chi_s p)=L\in [0,1)$, then
\begin{equation}
\label{eq:sublinear-dist}
\lim_{s\to \infty}\frac{d_{\Sigma_s}(\chi_{s}p,\chi_{s}q)}{s}=\sqrt{\frac{L}{1-L}}.
\end{equation}
In particular, $L\neq 0$ if and only if there exists $c>0$ such that $d_{\Sigma_s}(\chi_s p,\chi_s q)\geq cs$ for all sufficiently large $s$. $L=0$ if and only if $d_{\Sigma_s}(\chi_s p,\chi_s q)=o(s)$ for all sufficiently large $s$. 
\end{proposition}
\begin{proof}
Write $p_s=\chi_s p$, $q_s=\chi_s q$, so that $R(p_s)\to L\in[0,1)$ and $R(q_s)\to0$. 
Similar to \cite[Lemma~7.1]{HS26}, we have 
\begin{equation}
\label{eq:o-distance-limit}
\frac{d_g(o,p_s)}s\to\frac1{\sqrt{1-L}},
 \qquad
 \frac{d_g(o,q_s)}s\to1.
\end{equation}
Fix $\delta\in(0,1)$.
In what follows $C$ is a positive constant that may change from line to line, but stays independent of $\delta,t,s,b$. 
Let $\delta s\le t\le s$ and $\gamma_{s,t}:[0,d_g(p_s,q_t)]\to M$ be a unit-speed
minimizing geodesic from $\gamma_{s,t}(0)=q_t$ to $\gamma_{s,t}(d_g(p_s,q_t))=p_s$. The function $t\mapsto d_g(p_s,q_t)^2$ is locally Lipschitz and by the first variation formula, we obtain for almost every $t$, 
$$\frac{d}{dt}d_g(p_s,q_t)^2=-\frac{2d_g(p_s,q_t)}{|\nabla \rho|^2(q_t)} \langle \nabla \rho(q_t),\g_{s,t}'(0)\rangle.$$

\begin{claim}
\label{claim:zero-tip-sublinear}
We have as $s\to \infty$, 
\[
\sup_{\delta s\le t\le s}
\frac{|d_g(p_s,q_t) \langle \nabla \rho(q_t),\g_{s,t}'(0)\rangle-(s-t)|}{s}\to 0.
\]
\end{claim}
\begin{proof}
To simplify the notation, we write $\ell=d_g(p_s,q_t)$, $\gamma=\gamma_{s,t}$, $w=(\rho \circ \gamma)'$, so that $|w|\le 1$. Fix $b>0$. Suppose that $\ell>2b$. Then, \cite[Theorem~17.4]{Ham95} gives 
\[
 \int_b^{\ell-b}w'=\int_b^{\ell-b}\Ric(\g',\g')
 \le \frac{2(4-1)}b.
\]
The soliton identity gives $|\nabla R|\le2R$, hence $\int_0^b w'\,du\le b e^{2b}R(q_t)$. 
Since $w'\ge0$, for $u\in [0,\ell-b]$, we have
\[
 0\le w(u)-w(0)=\int_0^u w'
 \le \int_0^b w'+\int_b^{\ell-b}w'
 \le b e^{2b}R(q_t)+\frac{6}b.
\]
Since $|w|\le1$, we get $0\le w(u)-w(0)\le2$ on $[\ell-b,\ell]$. Integrating the above inequality on $[0,\ell]$, 
\[
\begin{aligned}
 0\le s-t-\ell w(0)
 &=\int_0^\ell \bigl(w(u)-w(0)\bigr)\,du\le\ell\left(b e^{2b}R(q_t)+\frac{6}b\right)+2b,
\end{aligned}
\]
where we used $\int_0^\ell w\,du
=\rho(p_s)-\rho(q_t)=s-t$. Since $\ell=O(s)$ and $R(q_t)\leq R(q_{\delta s})$, we obtain for all sufficiently large $s$, 
\[
0\le\frac{s-t-\ell w(0)}s
\le C\left(b e^{2b}R(q_{\delta s})+\frac1b\right)+\frac{2b}s.
\]
On the other hand, if $\ell\leq 2b$, then $0\leq s-t-\ell w(0)\leq 2\ell\leq 4b$. In both cases, using $R(q_{\delta s})\to 0$, we get for all sufficiently large $s$, 
\[
\sup_{\delta s\le t\le s}
\frac{|\ell w(0)-(s-t)|}{s}\le \frac{C}{b}.
\]
Since $b>0$ was arbitrary, this proves Claim~\ref{claim:zero-tip-sublinear}. 
\end{proof}
From Claim~\ref{claim:zero-tip-sublinear}, we get, as $s\to \infty$, 
$$\sup_{\delta s\leq t\le s}\left|\frac{|\nabla \rho|^2 (q_t)}{2s}\frac{d}{dt}d_g(p_s,q_t)^2+\frac{s-t}{s}\right|\to 0.$$
Since $|\nabla\rho(q_t)|^2=1-R(q_t)\to1$ uniformly for $\delta s\leq t\leq s$, we obtain 
$$\sup_{\delta s\leq t\le s}\left|\frac{1}{2s}\frac{d}{dt}d_g(p_s,q_t)^2+\frac{s-t}{s}\right|\to 0.$$
Integrating in $t$ from $\delta s$ to $s$, we obtain as $s\to \infty$, 
$$\left|\frac{d_g(p_s,q_{s})^2-d_g(p_s,q_{\delta s})^2}{s^2}+(1-\delta)^2\right|\to 0.$$
Since $|d_g(p_s,q_{\delta s})-d_g(o,p_s)|\le d_g(o,q_{\delta s})$, we may use \eqref{eq:o-distance-limit} to obtain as $s\to \infty$, 
\[
 \limsup_{s\to\infty}
 \left|
  \frac{d_g(p_s,q_s)^2}{s^2}
  -\frac1{1-L}+(1-\delta)^2
 \right|\le C\delta.
\]
Letting $\delta\to 0$ gives
\begin{equation}
\label{eq:d-g-p-s-q-s}
\lim_{s\to \infty}\frac{d_g(p_s,q_{s})^2}{s^2}=\frac{L}{1-L}.
\end{equation}
Let $\epsilon>0$. From Claim~\ref{claim:zero-tip-sublinear}, we also obtain as $s\to \infty$, 
$$\frac{|d_g(p_s,q_s) \langle \nabla \rho(q_s),\g_{s,s}'(0)\rangle|}{s}\to 0.$$
Since $u\mapsto \rho(\gamma_{s,s}(u))$ is convex, we have 
$$\rho(\gamma_{s,s}(u))\geq \rho(q_s)+(\rho\circ \g_{s,s})'(0)\cdot u\geq s-o(s).$$
For all sufficiently large $s$, it follows that $\rho(\gamma_{(1+\epsilon)s,(1+\epsilon)s})>s$. If $s_2\geq s_1$, then $\chi_{s_1}\circ\chi_{s_2}^{-1}:\Sigma_{s_2}\to \Sigma_{s_1}$ is a diffeomorphism such that 
\[
(\chi_{s_1}\circ\chi_{s_2}^{-1})^*
(g|_{T\Sigma_{s_1}})
\le g|_{T\Sigma_{s_2}},
\qquad s_2\ge s_1.
\]
See for instance, the proof of \cite[Proposition~2.2]{Deruelle12}. Let $\sigma_s:[0,d_g(p_{(1+\epsilon)s},q_{(1+\epsilon)s})]\to \Sigma_s$ denote the curve obtained by projecting the geodesic $\gamma_{(1+\epsilon)s,(1+\epsilon)s}$ onto $\Sigma_s$, i.e., $\sigma_s(t)=\chi_{s}\circ \chi_{\rho(\gamma_{(1+\epsilon)s,(1+\epsilon)s}(t))}^{-1}\gamma_{(1+\epsilon)s,(1+\epsilon)s}(t)$.  Then, $\sigma_s$ is a smooth curve lying in $\Sigma_s$ such that $\sigma_s(0)=q_s$, $\sigma_s(d_g(p_{(1+\epsilon)s},q_{(1+\epsilon)s}))=p_s$, and $|\sigma_s'|\leq 1$. Thus, for all sufficiently large $s$, 
$$d_g(p_s,q_s)\leq d_{\Sigma_s}(p_s,q_s)\leq \int_0^{d_g(p_{(1+\epsilon)s},q_{(1+\epsilon)s})}|\sigma_s'|\,dt\le d_g(p_{(1+\epsilon)s},q_{(1+\epsilon)s}).$$
Dividing by $s$, taking $s\to \infty$ and using \eqref{eq:d-g-p-s-q-s}, we obtain
$$\frac{L}{1-L}\leq \liminf_{s\to \infty}\frac{d_{\Sigma_s}(p_s,q_s)^2}{s^2}\leq \limsup_{s\to \infty}\frac{d_{\Sigma_s}(p_s,q_s)^2}{s^2}\leq \frac{L}{1-L}(1+\epsilon)^2.$$
Taking $\epsilon\to 0$ proves Proposition~\ref{prop:zero-tip-sublinear}. 
\end{proof}

\section{Dichotomy of modes}
\label{sec:dichotomy}
In this section, our goal is to set up the mode analysis for the radius of symmetry. We closely follow \cite{MMS23} and \cite[Section~3]{ABDS22} without using $R\to 0$ at infinity, but instead using properties derived in \cite{HS26}. For completeness and ease of reference, we present many of the arguments adapted from \cite{MMS23}, \cite{ABDS22}. 

\subsection{Reference point and bounds on errors}

\begin{proposition}
\label{prop:reference-bound}
After increasing $s_0$ if necessary, there exists $q\in\Sigma_{s_0}$ such that
\begin{equation}
\label{eq:ref-1}
\lim_{u\to\infty}uR(\Phi_{-u}(q))=1.
\end{equation}
\end{proposition}

\begin{proof}
Let $\epsilon_*,N_*$ and $\epsilon_1,r_1$ be the constants in \cite[Theorem~4.7(4)]{HS26} and \cite[Theorem~5.2]{HS26}, respectively. Choose $\epsilon<\min\{\epsilon_*,\epsilon_1\}$. 
By \cite[Lemma~2.8(iv)]{HS26}, after increasing $s_0$ we may
choose an $(\epsilon,2)$-center $q\in\Sigma_{s_0}$ with $d_g(o,q)>\max\{N_*,r_1\}.$
Theorem~4.7(4) of \cite{HS26} shows that the dimension reduction along $\Phi_{-u}q$ yields $\sph^2\times \R$. \cite[Theorem~5.2]{HS26} then gives the result.
\end{proof}

\begin{proposition}
\label{prop:reference-cylinder}
For this choice of $q$ and any sequence $t_j\to-\infty$,
\begin{equation}
\label{eq:ref-2}
(M,(-t_j)^{-1}g(t_j),q)\to 
 (\mathbb S^2\times\mathbb R^2,2g_{\mathbb S^2}+g_{\mathbb R^2}),
\end{equation}
in the smooth Cheeger--Gromov sense. The corresponding flows $(M,(-t_j)^{-1}g((-t_j)t),q)$ converge to $(\sph^2\times \R^2,-2tg_{\sph^2}+g_{\R^2})$ on compact time intervals in $t<0$. 
\end{proposition}
\begin{proof}
Put $u_j=-t_j$, $q_j=\Phi_{-u_j}(q)$, and $R_j=R(q_j)$. By \cite[Theorem 4.7(4)]{HS26}, every subsequence has a further subsequence for which
\[
(M,R_jg,q_j)\to 
(\mathbb S^2\times\mathbb R^2,
 2g_{\mathbb S^2}+g_{\mathbb R^2}).
\]
Proposition~\ref{prop:reference-bound} gives $u_jR_j\to1$. Since $(M,u_j^{-1}g(-u_j),q)$ is isometric to $(M,(u_jR_j)^{-1}R_jg,q_j),$
the claimed convergence follows. For the convergence of flows, apply the spacetime dimension reduction \cite[Proposition 4.2]{HS26} at $\Phi_{u_jb}(q)$ for $b<0$. 
\end{proof}

\begin{proposition}
\label{prop:level-reference}
For this choice of $q$, as $s\to\infty$,
\begin{equation}
\label{eq:ref-3}
 (M,s^{-1}g,\chi_s(q))\to 
 (\mathbb S^2\times\mathbb R^2,2g_{\mathbb S^2}+g_{\mathbb R^2},q_\infty),
\end{equation}
locally smoothly. 
\end{proposition}
\begin{proof}
Write $\chi_s(q)=\Phi_{t_s}(q)$. The preceding proposition gives $R\to0$, and therefore $|\nabla\rho|^2\to1$, as $t\to-\infty$. Integrating $\frac{d}{dt}\rho(\Phi_t(q))=-|\nabla\rho|^2$ gives $-t_s/s\to 1$. The statement now follows from Proposition~\ref{prop:reference-cylinder}. 
\end{proof}
Any point $q\in \Sigma$ satisfying \eqref{eq:ref-1}, \eqref{eq:ref-2}, and \eqref{eq:ref-3} will be called a reference point. Choose any reference point $q$, set $q_s=\chi_s(q)$, and let $\mathcal O_s=O(3)\cdot q_s\subset\Sigma_s$ be the sphere of symmetry passing through $q_s$, also called the spherical orbit of $q_s$. Set $q_s=\chi_s(q)$ and define
\begin{equation}\label{eq:tip-distances-definition}
 d_\pm(s):=d_{\Sigma_s}(p_s^\pm,q_s).
\end{equation}
By rotational symmetry, this equals the distance from $p_s^\pm$ to the spherical orbit through $q_s$. 

Using $O(3)$-symmetry, pulling the coordinates on $\Sigma_s$ back by $\chi_s$, we have on $\Sigma$, 
\begin{equation}\label{eq:metric-and-error}
 \bar g_s=dz^2+F(z,s)^2g_{\mathbb S^2},\qquad
 \partial_s\bar g_s=2\Ric_{\bar g_s}+2\bar E_s.
\end{equation}
Here $F(\cdot,s):\tilde I_s\to \R$ where $\tilde I_s:=(-d_-(s),d_+(s))$ and it extends continuously to the closure of ${\tilde I}_s$ with $F(-d_-(s),s)=F(d_+(s),s)=0.$ The arclength coordinate $z=z(x,s)$ depends on $s$. For scalar functions written in the variables $(z,s)$, a subscript $s$ denotes differentiation at fixed $z$. One can compute 
\begin{equation}
\label{eq:error-in-variation}
 \bar E_s=\chi_s^*\left(\frac{\Ric_g|_{T\Sigma_s}}{|\nabla\rho|^2}\right)
             -\Ric_{\bar g_s}.
\end{equation}
Choose a $g$-orthonormal frame with $e_3=\partial_z$, $e_4=n$, and $e_1,e_2$ tangent to the spherical orbits. Set
\[
 R_{ijkl}:=\Rm(e_i,e_j,e_k,e_l),\qquad
 R_{ij}:=\Ric_g(e_i,e_j),
\]
using the convention $R_{ijji}=\sec_g(e_i, e_j)$ for $i\neq j$. Bars denote intrinsic quantities on the level set. The Gauss equation gives 
\begin{equation}
\label{eq:gauss}
\bar R_{ij}
 =R_{ij}-R_{4ij4}
   +\frac{R_{ij}(R-R_{44})-\sum_{k=1}^3R_{ik}R_{kj}}{1-R}.
\end{equation}
We will write $ \bar E_s=E^{\rm rad}dz^2+E^{\rm orb}g_{\mathbb S^2}$. For tangential indices $i,j\in\{1,2,3\}$,  write $E_{ij}$ for the components of $\bar E_s$ in the frame $\{e_1,e_2,e_3\}$. In particular,
\[
 E^{\rm rad}=E_{33},\qquad
 E^{\rm orb}=F^2E_{11}=F^2E_{22}.
\]
From \eqref{eq:error-in-variation} and the Gauss equation, it follows that 
\begin{equation}
\label{eq:error-identity}
 E_{ij}=R_{4ij4}
   +\frac{R_{ij}R_{44}+\sum_{k=1}^3R_{ik}R_{kj}}{1-R}.
    \end{equation}
Recall the definition of $(\epsilon,2)$-centers from \cite{HS26}. 
\begin{lemma}
\label{lem:large-orbit-convergence}
On $M$, wherever $F>0$, we have $ R\le C F^{-2}.$ 
For every $a>0$ and $\epsilon>0$, there exists $s_*(a,\epsilon)$ such that every point in $\{\rho=s,\ F\geq a\sqrt{s}\}$ is an $(\epsilon,2)$-center whenever $s\geq s_*(a,\epsilon)$. 
\end{lemma}
\begin{proof}
The estimate $RF^2\leq C$ follows from \cite[Lemma~2.4]{MMS23}, whose proof uses only dimension reduction. Suppose the second statement of the lemma fails. 
Then for some $a,\epsilon>0$ there are $s_j\to\infty$ and $x_j\in\Sigma_{s_j}$ such that $F(x_j)\geq a\sqrt{s_j},$ 
but $x_j$ is not an $(\epsilon,2)$-center. 
In particular, after passing to a subsequence $(M,R(x_j)g,x_j)$ converges to $\Bry^3\times \R$. 
Moreover, $R(x_j)\leq CF(x_j)^{-2}\leq Ca^{-2} s_j^{-1}$ showing that $R(x_j)\to 0$. 
In particular, $(\Sigma_{s_j},R(x_j)g,x_j)$ converges to the Bryant soliton $(\Bry^3,\tilde{g}_0,x_\infty)$. 
By smooth convergence, since $\Sigma_{s_j}$ has at most two tips for all large $j$, it follows that for some $p_j\in \{p_{s_j}^+,p_{s_j}^-\}$, we have $\sup_j R(x_j)d_g(x_j,p_j)^2<\infty$. By Perelman's long-range estimates, $R(p_j)\leq CR(x_j)\leq Ca^{-2} s_j^{-1}$. This is a contradiction to Lemma~\ref{lem:poles}. 
\end{proof}

\begin{lemma}
\label{lem:error-sign}
$E^{\rm rad}\ge0$ and $E^{\rm orb}\ge0$ on every sufficiently large level.
\end{lemma}
\begin{proof}
From \eqref{eq:error-identity}, we obtain 
\begin{equation}
\label{eq:Erad-Etan-exp}
{
\begin{aligned}
E^{\rm rad}
=R_{4334}+\frac{R_{33}(R_{44}+R_{33})}{1-R},\qquad E^{\rm orb}
=F^2\left(
R_{4114}+\frac{R_{11}(R_{44}+R_{11})}{1-R}
\right).
\end{aligned}}
\end{equation}
The lemma now follows. 
\end{proof}

\begin{lemma}
\label{lem:warping-derivatives}
For sufficiently large $s$, wherever $F>0$,
\[
 |F_z|\leq1,\qquad F_{zz}<0,\qquad
 F^m|\partial_z^{m+1}F|\leq C_m,
\]
for every integer $m\geq 0$.
\end{lemma}

\begin{proof}
Since $\sec_g\ge0$, we get $\sec_{\bar g_s}>0$ using Gauss equations. Moreover, 
\[
{
\sec_{\bar g_s}(e_1,e_2)=\frac{1-F_z^2}{F^2},
\qquad
\sec_{\bar g_s}(\partial_z,e_1)=-\frac{F_{zz}}F.
}
\]
This proves the first two assertions. Using dimension reduction, we obtain
\[
 |\bar\nabla^k\overline{\Rm}|\le C_kR^{1+k/2}
                   \le C_kF^{-2-k}.
\]
In particular $F|F_{zz}|\le C$. Differentiating the identity $\bar R=2F^{-2}(1-F_z^2-2FF_{zz})$ and using an induction argument on $m$ proves the lemma.  
\end{proof}

\begin{lemma}
\label{lem:error-size}
For all sufficiently large $s$, wherever $F>0$, we have 
\[
 |E^{\rm rad}|\le CF^{-4},\qquad
 |\partial_z^kE^{\rm orb}|\le C_kF^{-2-k},
\]
for every integer $k\geq 0$.
\end{lemma}
\begin{proof}
By Lemma~\ref{lem:large-orbit-convergence}, we have $|\Ric|\leq R\leq CF^{-2}$ and recall that $|\nabla^k\Rm|\le C_kR^{1+k/2}$ on $M$, whenever $F>0$. Recall the soliton identity in \cite[(3.11)]{MMS23}, also \cite[Lemma~3.1]{DZ20}, for $i,j\in\{1,2,3\}$, 
\begin{equation}\label{eq:normal-curvature}
 (1-R)R_{4ij4}
 =-\frac12(\nabla^2R)_{ij}
  -\sum_{k=1}^4R_{ik}R_{kj}
  +(\Delta\Ric)_{ij}
  +2\sum_{k,l=1}^4R_{iklj}R_{kl}.
\end{equation}
Combining this with \eqref{eq:Erad-Etan-exp} and using the uniform lower bound for $1-R$ away from $o$, we obtain $|E^{\rm rad}|\le CF^{-4}$.  Differentiating \eqref{eq:normal-curvature} and using  \eqref{eq:Erad-Etan-exp}, we obtain $|\bar\nabla^k\bar E_s|\le C_kR^{2+k/2}\le C_kF^{-4-k}.$
Recall that $\bar \nabla_{\partial_z}\partial_z=0$, $\bar \nabla_{\partial_z}e_1=0$, and $\partial_z^k\!\left(F^{-2}E^{\rm orb}\right)
=(\bar\nabla^k\bar E_s)
({\partial_z,\ldots,\partial_z},e_1,e_1).$
Lemma~\ref{lem:warping-derivatives} gives $|\partial_z^j(F^2)|\le C_jF^{2-j}.$ Hence 
\[
|\partial_z^kE^{\rm orb}|
\le\sum_{j=0}^k\binom kj
C_jF^{2-j}C_{k-j}F^{-4-(k-j)}
\le C_kF^{-2-k}.
\]
\end{proof}

\subsection{Estimates on the rescaling of the radius function}

\begin{proposition}[{\cite[Proposition 3.16]{MMS23}}]\label{prop:F-equation}
In coordinates $(z,s)$, 
\begin{align}\label{eq:F-equation}
 -F_s={}&F_{zz}-F^{-1}(1+F_z^2)
 +2F_z\left\{\frac{F_z}{F}(0,s)
                -\int_0^z\frac{F_z^2}{F^2}(\zeta,s)\,d\zeta\right\}
 -F^{-1}E^{\rm orb}+F_z\int_0^zE^{\rm rad}(\zeta,s)\,d\zeta.
\end{align}

\end{proposition}
\begin{proof}
Differentiate \eqref{eq:metric-and-error} and use \eqref{eq:Erad-Etan-exp}. The calculation of \cite[Proposition 3.16]{MMS23} goes unchanged. 
\end{proof}

\begin{lemma}
\label{lem:rmax}
For $r_{\max}(s):=\max_{z\in \tilde I_s} F(z,s)$, 
\[
 2s-C\le r_{\max}(s)^2\le 2s(1+o(1)).
\]
\end{lemma}
\begin{proof}
Let $z_*(s)$ be the unique maximum point of $F(\cdot,s)$. 
At this point, $F_z=0$, and \eqref{eq:F-equation} gives 
\[
 \frac{d}{ds}r_{\max}(s)^2
 =2-2FF_{zz}+2E^{\rm orb}\geq2.
\]
Integration yields $r_{\max}(s)^2\geq2s-C$. As a consequence of Lemma~\ref{lem:large-orbit-convergence}, smooth convergence at $z_*(s)$ implies $FF_{zz}=o(1)$ there, while Lemma~\ref{lem:error-size} gives $E^{\rm orb}=o(1)$. Hence $\frac{d}{ds}r_{\max}(s)^2=2+o(1).$ Integrating again proves the upper bound. 
\end{proof}
For universal constants $a\gg 1$ and $r_*\in(0,1)$, let $\psi_a$ denote the barrier functions of \cite[Section 2]{Bre20}, as recorded in \cite[Proposition~2.19]{MMS23}, which are real-valued $C^1$ functions on $[r_*a^{-1},1+\frac{a^{-2}}{100}]$ that satisfy a certain ordinary differential inequality. We use, in particular, $ \psi_a(u)\le Ca^{-2}$ for $u\in [\frac{1}{10},1+\frac{a^{-2}}{100}]$ and $ \psi_a(r_*a^{-1})\ge3/2.$ 

\begin{proposition}[{\cite[Proposition 3.19]{MMS23}}]\label{prop:barrier}
There exist constants $a_0,D>0$ with the following property. Suppose $a\geq a_0$ and $S\geq Da^3$, and assume that
\[
 \frac{r_{\max}(s)}{\sqrt{2s}}
 \leq1+\frac{a^{-2}}{100}
 \qquad\text{for all }s\geq S.
\]
Then
\[
 F_z(z,s)^2
 \leq\psi_a\left(\frac{F(z,s)}{\sqrt{2s}}\right),
\]
whenever $s\geq S$ and $F(z,s)\geq r_*a^{-1}\sqrt{2s}$. 
\end{proposition}
\begin{proof}
The argument of \cite[Proposition 3.19]{MMS23} goes through and doesn't require $R\to 0$ at infinity. 
\end{proof}

Writing $\xi:=z/\sqrt{s}$, the equation $\bar g_s=dz^2+F(z,s)^2g_{\mathbb S^2}$ becomes, in the $\xi$-coordinate, 
\[
{
s^{-1}\bar g_s
=d\xi^2+
\left(\frac{F(\sqrt s\,\xi,s)}{\sqrt s}\right)^2g_{\mathbb S^2}.
}
\]
In the limit, when $|\xi|$ is bounded, this expression becomes $d\xi^2+2g_{\mathbb S^2}$. This motivates the use of the quantity $(\xi,s)\mapsto \frac{F(\sqrt s\,\xi,s)}{\sqrt s}-\sqrt{2}$. 

Define, as in \cite{MMS23}, 
\begin{gather}
 s=e^{-\tau},\quad z=e^{-\tau/2}\xi,\quad
 G(\xi,\tau)=e^{\tau/2}F(e^{-\tau/2}\xi,e^{-\tau})-\sqrt2,
       \label{eq:G-definition}\\
 \mathcal E^{\rm rad}(\xi,\tau)=e^{-\tau}E^{\rm rad}(z,s),\qquad
 \mathcal E^{\rm orb}(\xi,\tau)=E^{\rm orb}(z,s).
       \label{eq:rescaled-errors}
\end{gather}
Here $-\infty<\tau\leq -\log s_0$ and $G(\cdot,\tau):J_\tau\to \R$ where $J_\tau:=(-e^{\tau/2}d_-(e^{-\tau}),e^{\tau/2}d_+(e^{-\tau}))$. 
$G(\cdot,\tau)$ continuously extends to $\bar J_\tau$ with $G(-e^{\tau/2}d_-(e^{-\tau}),\tau)=G(e^{\tau/2}d_+(e^{-\tau}),\tau)=-\sqrt2.$
We have for all sufficiently negative $\tau$, 
$$-\sqrt{2}\leq G(\cdot,\tau)\leq e^{\tau/2}r_{\max}(e^{-\tau})-\sqrt{2}=o(1).$$
Moreover, $G\to0$ locally smoothly as $\tau\to-\infty$. 

\begin{proposition}[{\cite[Proposition 3.24]{MMS23}}]\label{prop:G-equation}
We have 
\begin{align}\label{eq:G-equation}
 G_\tau={}&G_{\xi\xi}-\frac\xi2G_\xi
 +\frac12(\sqrt2+G)-(\sqrt2+G)^{-1}(1+G_\xi^2)\notag\\
 &+2G_\xi\left\{\frac{G_\xi(0,\tau)}{\sqrt2+G(0,\tau)}
    -\int_0^\xi\frac{G_\xi(x,\tau)^2}{(\sqrt2+G(x,\tau))^2}\,dx\right\}
 -(\sqrt2+G)^{-1}\mathcal E^{\rm orb}
 +G_\xi\int_0^\xi\mathcal E^{\rm rad}(x,\tau)\,dx.
\end{align}
\end{proposition}
\begin{proof}
This follows by rescaling \eqref{eq:F-equation}. 
\end{proof}
Using the identity $\frac{\sqrt2+G}{2}-\frac1{\sqrt2+G}
=G-\frac{G^2}{2(\sqrt2+G)}$, Proposition~\ref{prop:G-equation} shows that the evolution of $G(\cdot,\tau)$ has the form 
\[
G_\tau
={G_{\xi\xi}-\frac{\xi}{2}G_\xi+G}
+Q,
\]
where $Q=Q_{1a}+Q_{1b}+Q_2$ consists of nonlinear terms and error terms: 
\begin{equation}
\label{eq:decomposition-of-errors}
\begin{aligned}
Q_{1a}={}&-\frac{G^2}{2(\sqrt{2}+G)}-\frac{G_\xi^2}{(\sqrt{2}+G)}
+\frac{2G_\xi G_\xi(0)}{\sqrt{2}+G(0)},\\
Q_{1b}={}&-2G_\xi \int_0^\xi\frac{G_\xi(x)^2}{(\sqrt{2}+G(x))^2}\,dx,\\
Q_2={}&-\frac{\mathcal E^{\rm orb}}{\sqrt{2}+G}
+G_\xi\int_0^\xi\mathcal E^{\rm rad}(x)\,dx.
\end{aligned}
\end{equation}
We define the linear differential operator 
$$\mathcal L:=\partial_{\xi\xi}-\frac\xi2\partial_\xi+1,$$
so that $Q=(\partial_\tau-\mathcal L)G$. Set
\[
 \varrho_{\max}(\tau)=e^{\tau/4}+\max_{\xi\in J_\tau} G(\xi,\tau),\qquad
 \varrho(\tau)=\sup_{\sigma\le\tau}\varrho_{\max}(\sigma),\qquad
 \delta(\tau)=\varrho(\tau)+\sup_{\sigma\le\tau}|G(0,\sigma)|.
\]
(To avoid conflict with the potential $\rho$, write $\varrho$ for the
quantity denoted by $\rho$ in \cite[Definition 3.25]{MMS23}. ) For sufficiently negative $\tau$, by Lemma~\ref{lem:rmax}, 
\[
 \varrho_{\max}(\tau)= e^{\tau/4}+e^{\tau/2}r_{\max}(e^{-\tau})-\sqrt{2}\ge e^{\tau/4}+ \sqrt{2-Ce^\tau}-\sqrt{2}\ge \tfrac12e^{\tau/4},
\]
and $\delta(\tau)\ge\varrho(\tau)\ge\tfrac12e^{\tau/4}.$ Moreover, Lemma~\ref{lem:rmax} also shows that both $\varrho(\tau),\delta(\tau)\to0$ uniformly as $\tau\to -\infty$. 

\begin{proposition}
\label{prop:rescaled-barrier}
For any sufficiently negative $\tau_*$, the following property holds true. Set $a=10^{-1}\varrho(\tau_*)^{-1/2}$. Then for all $\tau\leq \tau_*$, 
\[
 G_\xi^2\le\psi_a(1+G/\sqrt2)
 \quad\hbox{on }\{\sqrt2+G\ge\sqrt2r_*a^{-1}\}.
\]
On the smaller region $\{1+G/\sqrt2\ge1/10\}$ and for all $\tau\leq \tau_*$, one also has
$|G_\xi|\le C\sqrt{\varrho(\tau_*)}.$
Here $C$ is universal.  
\end{proposition}

\begin{proof}
Set $a=10^{-1}\varrho(\tau_*)^{-1/2}, S=e^{-\tau_*}$, for $\tau_*$ which we will now choose. The lower bound $\varrho(\tau_*)\geq\frac12e^{\tau_*/4}$ gives $a^3\leq Ce^{-3\tau_*/8}$. Choose any $\tau_*\ll 0$ to ensure $S\geq Da^3$ and $a\geq \max(a_0,10)$. For every $\sigma\leq\tau_*$,
\[
 \frac{r_{\max}(e^{-\sigma})}{\sqrt{2e^{-\sigma}}}
 =1+\frac{\max_\xi G(\xi,\sigma)}{\sqrt2}
 \leq1+\frac{\varrho(\tau_*)}{\sqrt2}
 \leq1+\frac{a^{-2}}{100}.
\]
Proposition~\ref{prop:barrier}, applied on $[S,\infty)$, therefore proves the first assertion. On the smaller region $1+G/\sqrt2\geq1/10$, the argument of $\psi_a$ belongs to $[1/10,1+a^{-2}/100]$. Thus $\psi_a(1+G/\sqrt2)\leq Ca^{-2}\leq C\varrho(\tau_*)$. 
\end{proof}

\begin{lemma}
\label{lem:central-domain}
Let $I_\tau=[-\delta(\tau)^{-1/100},\delta(\tau)^{-1/100}]$. For sufficiently negative $\tau$, the function $\xi\mapsto G(\xi,\tau)$ is defined on $\xi\in I_\tau$, and there $G\ge-1.$ In particular, $1+ G(\cdot,\tau)/{\sqrt2}>1/10$ on $I_\tau$. 
\end{lemma}

\begin{proof}
Since $G(0,\tau)\to0$, fix any $\tau$ such that $G(0,\tau)\ge-\frac14$ and such that Proposition~\ref{prop:rescaled-barrier} can be applied with $\tau_*:=\tau$. 
We claim that $G(\xi,\tau)\ge-1$ for all $0\leq \xi\leq \varrho(\tau)^{-1/100}$. Suppose not and let $\xi_1\le \varrho(\tau)^{-1/100}$  be the first point where $G(\xi_1,\tau)=-1$. Then, $G(\xi,\tau)\ge-1$ for $0\leq \xi\le \xi_1$. Then, for all $0\leq \xi\leq \xi_1$, 
\[
1+\frac{G(\xi,\tau)}{\sqrt2}
\ge1-\frac1{\sqrt2}>\frac1{10}.
\]
Proposition~\ref{prop:rescaled-barrier} shows that $|G_\xi|\le C\sqrt{\varrho(\tau)}$  for all $0\leq \xi\leq \xi_1$. Integrating, we get 
\[
\begin{aligned}
G(\xi_1,\tau)=G(0,\tau)+\int_0^{\xi_1}G_\xi(x,\tau)\,dx\ge-\frac14-C\sqrt{\varrho(\tau)}\,\xi_1\ge -\frac14-C{\varrho(\tau)}^{49/100}.
\end{aligned}
\]
If $\tau$ is sufficiently negative, the right hand side is atleast $-1/2$ contradicting the definition of $\xi_1$. The same argument applies on $[-\varrho(\tau)^{-1/100},0]$, proving the lemma. 
\end{proof}
Recall that by definition of $G$, we have $|G_\xi(\xi,\tau)|\le 1$ on $J_\tau$. 
\begin{lemma}
\label{lem:higher-derivatives}
For each $m\geq 0$, for sufficiently negative $\tau$, the following estimates hold true on $I_\tau$:  
$$\begin{aligned}
|G|+|G_\xi|&\le C\delta(\tau)^{1/4}\\
|G_{\xi\xi}|&\le C\delta(\tau)^{1/8}\\
|\partial_\xi^mG|&\le C_m.
\end{aligned}$$
Here $C,C_m$ are independent of $\tau$ once $\tau\ll 0$. 
\end{lemma}
\begin{proof}
Combining Lemma~\ref{lem:central-domain} with Proposition~\ref{prop:rescaled-barrier}, we obtain the $C^1$-bound $|G_\xi|\le C\sqrt{\varrho}\le C\delta^{1/2}.$ Using $|G(0,\tau)|\le\delta$, we obtain the $C^0$-bound: 
\[
|G(\xi,\tau)|
\le\delta+C\delta^{1/2}|\xi|
\le\delta+C\delta^{49/100}
\le C'\delta^{1/4}.
\]
For the $C^m$-bounds, rescale the estimates in Lemma~\ref{lem:warping-derivatives}. Using $\partial_\xi^mG=(\sqrt{s})^{m-1}\partial_z^mF$ and the estimates $F^{m-1}|\partial_z^mF|\le C_{m-1},$ we obtain on $J_\tau$,  
\begin{equation}
\label{eq:C-m-estimate-J-tau}
|\partial_\xi^mG|
\le C_{m-1}\left(\frac{\sqrt s}{F}\right)^{m-1}
=C_{m-1}(\sqrt{2}+G)^{-(m-1)}.
\end{equation}
This gives $C^m$-bounds on $I_\tau$ for all $m\geq 1$. At any $\xi\in I_\tau$, choose a direction pointing into the interval and Taylor’s formula for $h=\delta^{1/8}$ gives 
\[
|G_{\xi\xi}(\xi)|
\le\frac{2\|G_\xi\|_\infty}{h}
+\frac h2\|G_{\xi\xi\xi}\|_\infty
\le C\left(\frac{\delta^{1/4}}{\delta^{1/8}}
+\delta^{1/8}\right).
\]
Therefore $|G_{\xi\xi}|\le C\delta^{1/8}.$ 
\end{proof}
Put $d\nu=(4\pi)^{-1/2}e^{-\xi^2/4}\,d\xi$. 
\begin{lemma}
\label{lem:quartic}
For sufficiently negative $\tau$,
\begin{align*}
 |G_\xi(0,\tau)|^4
 &\le C\delta^{1/100}\int_{I_\tau}G^2\,d\nu,\\
 \int_{I_\tau}G_\xi^4\,d\nu
 &\le C\delta^{1/100}\int_{I_\tau}G^2\,d\nu
       +C\exp(-c\delta^{-1/50}).
\end{align*}
Here $C$ is independent of $\tau$ once $\tau\ll 0$. 
\end{lemma}

\begin{proof}
Let $\tau$ be sufficiently negative and write $a=|G_\xi(0,\tau)|$. The case $a=0$ is immediate, so we assume $a>0$. Lemma~\ref{lem:higher-derivatives} gives
\[
 G(x,\tau)-G(-x,\tau)=2G_\xi(0,\tau)x+O(x^3).
\]
In particular, $|G(x)-G(-x)|\ge ax$ for $0\le x\le c\sqrt a$, for some $c>0$. It follows that
\[
 \int_{I_\tau}G^2\,d\nu\ge c\int_0^{c\sqrt a}a^2x^2\,dx\ge c'a^{7/2}.
\]
Combining this with $|G_\xi|\le C\delta^{1/4}$, we get $a^4\le C\delta^{1/8}\int_{I_\tau}G^2\,d\nu$, showing the first estimate. For the second estimate, set $H=\int_{I_\tau}G^2\,d\nu$, $J=\int_{I_\tau} G_\xi^4\,d\nu$, $r=\delta(\tau)^{-1/100}$, and write $d\nu=w\,d\xi$. Integration by parts gives 
\[
 J=[GG_\xi^3w]_{-r}^r
 +\int_{-r}^r GG_\xi^2
       \left(\frac{\xi}{2}G_\xi-3G_{\xi\xi}\right)\,d\nu.
\]
The boundary term is $O(e^{-r^2/4})$, and the coefficient in parentheses is bounded by $C(r\delta^{1/4}+\delta^{1/8})\le C\delta^{1/8}$. Hence Cauchy--Schwarz and Young's inequality give
\[
 J\le Ce^{-r^2/4}+C\delta^{1/8}H^{1/2}J^{1/2}
 \le \tfrac12J+C\delta^{1/4}H+Ce^{-r^2/4}.
\]
This completes the proof. 
\end{proof}

\begin{lemma}
\label{lem:scaled-error}
Where $F>0$,
\[
 |\mathcal E^{\rm rad}|\le Ce^\tau(\sqrt2+G)^{-4},\qquad
 |\partial_\xi^k\mathcal E^{\rm orb}|
       \le C_ke^\tau(\sqrt2+G)^{-2-k}.
\]
In particular, all these quantities are $O_k(e^\tau)$ on $I_\tau$. 
\end{lemma}
\begin{proof}
The two estimates follow by rescaling the estimates in Lemma~\ref{lem:error-size}. On $I_\tau$, $\sqrt{2}+G\geq \sqrt{2}-1$, showing the final sentence. 
\end{proof}

\begin{lemma}
\label{lem:residual}
We have for sufficiently negative $\tau$, 
\[
 \int_{I_\tau}|Q|^2\,d\nu
 \le C\delta^{1/100}\int_{I_\tau}G^2\,d\nu
       +C\exp(-c\delta^{-1/50})+Ce^{2\tau}.
\]
\end{lemma}
\begin{proof}
Recall the decomposition $Q=Q_{1a}+Q_{1b}+Q_2$ in \eqref{eq:decomposition-of-errors}. Using the estimates on $G$, it follows that 
\[
\begin{aligned}
\int_{I_\tau}|Q_{1a}|^2\,d\nu
&\le C\left(\delta^{1/2}\int_{I_\tau}G^2\,d\nu+\int_{I_\tau}G^4_\xi\,d\nu+|G_\xi(0)|^4\right).
\end{aligned}
\]
Using the bound on the error terms, 
\[
\int_{I_\tau}|Q_2|^2\,d\nu\le Ce^{2\tau}.
\]
The bound on $Q_{1b}$ now follows by applying Lemma~\ref{lem:primitive-sobolev} to $h(\xi)=\int_0^\xi
\frac{G_\xi(x)^2}{(\sqrt2+G(x))^2}\,dx$: 
\begin{equation}
\label{eq:bound-Q-1b}
\int_{I_\tau}|Q_{1b}|^2\,d\nu\le C\delta^{1/2}\int_{I_\tau}G_\xi^4\,d\nu
\end{equation}
Combining these estimates with Lemma~\ref{lem:quartic} proves the lemma. 
\end{proof}

\subsection{Localization and dichotomy of modes}
\label{subsec:local-dichotomy}
Let $\mathcal H=L^2(\mathbb R,d\nu)$, where $d\nu=(4\pi)^{-1/2}e^{-\xi^2/4}\,d\xi$. The Hermite polynomials form a complete orthonormal eigenbasis $\{h_k\}_{k\ge0}$ for $\mathcal L$, with $\mathcal Lh_k=(1-k/2)h_k$ and
\[
 h_0=1,\quad h_1=\xi/\sqrt2,\quad
 h_2=(\xi^2-2)/(2\sqrt2),\quad
 h_3=(\xi^3-6\xi)/(4\sqrt3).
\]
Write $\mathcal H=\mathcal H_+\oplus\mathcal H_0\oplus\mathcal H_-$, where $\mathcal H_+=\operatorname{span}\{h_0,h_1\}$, $\mathcal H_0=\operatorname{span}\{h_2\}$, and $\mathcal H_-=(\mathcal H_+\oplus\mathcal H_0)^\perp$, and denote the corresponding orthogonal projections by $P_+,P_0,P_-$.

Choose an even function $\eta\in C_c^\infty(\mathbb R)$, $0\le\eta\le1$, equal to one on $[-1/2,1/2]$, supported in $[-1,1]$, and satisfying $x\eta'(x)\le0$. Define
\[
 \vartheta(\xi,\tau)=\eta(\delta(\tau)^{1/100}\xi),\quad
 \widehat G=\vartheta G,\quad
 \gamma=\|\widehat G\|_{\mathcal H}^2,\quad
 \gamma^\bullet=\|P_\bullet\widehat G\|_{\mathcal H}^2
 \quad(\bullet=+,0,-).
\]

\begin{lemma}
\label{lem:cutoff}
The function $\delta$ is locally Lipschitz and $0\le\delta'\le C$ almost everywhere for sufficiently negative $\tau$. 
We have $\partial_\tau\widehat G=\mathcal L\widehat G+\widehat Q,$ where, almost everywhere in $\tau$, 
\begin{equation}\label{eq:localized-residual}
 \|\widehat Q\|_{\mathcal H}^2
 \le C\delta^{1/100}\gamma+C\exp(-c\delta^{-1/50})+Ce^{2\tau}.
\end{equation}
\end{lemma}
\begin{proof}
At $z=0$, \eqref{eq:F-equation} gives
\[
 \frac{d}{ds}F(0,s)^2
 =2\bigl(1-F_z(0,s)^2-F(0,s)F_{zz}(0,s)
          +E^{\rm orb}(0,s)\bigr)
 =2+o(1),
\]
Combining this with Lemma~\ref{lem:large-orbit-convergence} 
\[
 \frac{d}{d\tau}\bigl(e^\tau F(0,e^{-\tau})^2\bigr)
 =e^\tau F(0,e^{-\tau})^2-\left.\frac{d}{ds}\right|_{s=e^{-\tau}}F(0,s)^2=o(1).
\]
This shows that $G(0,\cdot)$ and $\max_\xi G(\xi,\cdot)$ have uniformly bounded derivatives for sufficiently negative $\tau$. 
Hence $|G(0,\cdot)|$ and $\varrho_{\max}$ are uniformly Lipschitz. As a consequence, we have $0\le\delta'\le C$ almost everywhere. 

Differentiating $\widehat G=\vartheta G$ gives $\widehat Q=\vartheta Q+R_\vartheta,$ where
\[
R_\vartheta
=\vartheta_\tau G-2\vartheta_\xi G_\xi
-\vartheta_{\xi\xi}G+\frac{\xi}{2}\vartheta_\xi G.
\]
$\vartheta_\xi,\vartheta_\tau$ are both supported on $[-\delta(\tau)^{-1/100},-\delta(\tau)^{-1/100}/2]\cup [\delta(\tau)^{-1/100}/2,\delta(\tau)^{-1/100}]$ and 
\[
|\vartheta_\xi|\le C\delta(\tau)^{1/100},\qquad
|\vartheta_{\xi\xi}|\le C\delta(\tau)^{2/100},\qquad
|\xi\vartheta_\xi|\le C,\qquad
|\vartheta_\tau|\le C\delta(\tau)^{-1}.
\]
It follows that $|R_\vartheta|\le C\delta(\tau)^{-1}$ on $[\delta(\tau)^{-1/100}/2,\delta(\tau)^{-1/100}]$. Therefore,
\[
\begin{aligned}
\|R_\vartheta\|_{\mathcal H}^2\le C\delta(\tau)^{-2}\int_{|\xi|\ge \delta(\tau)^{-1/100}/2}d\nu\le C\delta(\tau)^{-2}e^{-c_0\delta(\tau)^{-1/50}}\leq C'e^{-c'\delta(\tau)^{-1/50}}
\end{aligned}
\]
Note that 
\[
\int_{I_\tau}G^2\,d\nu \leq \gamma+\int_{I_\tau}(1-\vartheta^2)G^2\,d\nu\leq \gamma+C\int_{|\xi|\ge \delta(\tau)^{-1/100}/2}d\nu\leq \gamma+Ce^{-c_0\delta(\tau)^{-1/50}}.
\]
Applying Lemma~\ref{lem:residual}, we obtain 
\[
 \|\vartheta Q\|_{\mathcal H}^2\leq \int_{I_\tau}|Q|^2\,d\nu
 \le C\delta^{1/100}\gamma+C\exp(-c\delta^{-1/50})+Ce^{2\tau}.
\]
Combining these estimates concludes the proof. 
\end{proof}
Set $p=8-\frac1{200}, b=\frac1{400}$ and 
\[
 S(\tau)=\sup_{u\in[\tau-1,\tau]}\gamma(u),\qquad
 T(\tau)=C\exp(-c\delta(\tau)^{-1/50})
                 +C\delta(\tau)^{-b}e^{2\tau}.
\]

\begin{lemma}
\label{lem:spectral-step}
For sufficiently negative $\tau$,
\begin{align*}
 \gamma^+(\tau-1)&\le e^{-1}\gamma^+(\tau)+C\delta(\tau)^bS(\tau)+T(\tau),\\
 |\gamma^0(\tau-1)-\gamma^0(\tau)|
     &\le C\delta(\tau)^bS(\tau)+T(\tau),\\
 \gamma^-(\tau-1)&\ge e\gamma^-(\tau)-C\delta(\tau)^bS(\tau)-T(\tau).
\end{align*}
\end{lemma}
\begin{proof}
Using Lemma~\ref{lem:cutoff}, we have almost everywhere, 
\[
\begin{aligned}
(\gamma^+)'(u)
&=2\langle P_+\widehat G,\mathcal LP_+\widehat G\rangle
+2\langle P_+\widehat G,P_+\widehat Q\rangle\ge\gamma^+(u)
-2\sqrt{\gamma^+(u)}\,\|\widehat Q(u)\|_{\mathcal H},
\end{aligned}
\]
because the eigenvalues on $\mathcal H_+$ are $1$ and $1/2$. Now fix $\tau$ and let $u\in[\tau-1,\tau]$. Since $\gamma^+(u)\le S(\tau)$, we obtain
\[
2\sqrt{\gamma^+(u)}\,\|\widehat Q(u)\|_{\mathcal H}
\le C\delta(\tau)^{1/200}S(\tau)
+C\bigl(e^{-c\delta(\tau)^{-1/50}}+e^\tau\bigr)\sqrt{S(\tau)}.
\]
By Young’s inequality, the right hand side is bounded above by $C\delta(\tau)^bS(\tau)+T(\tau).$ Thus $$(\gamma^+)'(u)\ge\gamma^+(u)-(C\delta(\tau)^bS(\tau)+T(\tau))$$ almost everywhere on $[\tau-1,\tau]$. Integrating gives  \[
\gamma^+(\tau)
\ge e\,\gamma^+(\tau-1)-(e-1)(C\delta(\tau)^bS(\tau)+T(\tau)),
\]
giving the first inequality. The other two inequalities follow similarly. 
\end{proof}

\begin{lemma}
\label{lem:max-second}
If $\xi_\tau^*$ is the unique maximum point of $G(\cdot,\tau)$, then $0<-G_{\xi\xi}(\xi_\tau^*,\tau)\le C\gamma(\tau)^{1/4}$. 
\end{lemma}
\begin{proof}
Suppose the domain of $G(\cdot,\tau)$ is $[-D_-,D_+]$. Since $|G_\xi|\leq 1$, we have $G(\xi^*_\tau)+\sqrt{2}\leq \min(D_+-\xi^*_\tau,\xi^*_\tau+D_-)$. Since $G(\xi_\tau^*,\tau)\to 0$, the interval $[\xi_\tau^*-1/2,\xi_\tau^*+1/2]$ lies within $J_\tau$. 

Fix a sufficiently negative $\tau$, and write $\xi_*=\xi_\tau^*$ and $a=-G_{\xi\xi}(\xi_*)>0$. Without loss of generality, we may assume $\xi_*\ge0$. By \eqref{eq:C-m-estimate-J-tau}, there exists $C_3>a$ such that $|G_{\xi\xi\xi}|\le C_3$ on $[\xi_*-1/2,\xi_*+1/2]$. 
Using Taylor's formula for $h=a/2C_3$,  we obtain 
\[
 G_\xi(\xi_*-h)\geq -G_{\xi\xi}(\xi_*)h-\frac{h^2}{2}C_3
 \ge\frac{a^2}{4C_3}.
\]
Using concavity, $ G_\xi(\xi)\ge\frac{a^2}{4C_3}$ for all $\xi\leq \xi_*-h$ and in particular on $J=[-2,-1]$. For $x,y\in J$, we have 
\[
 \int_J G^2\,d\xi
 \ge\frac12\int_J\int_J(G(x)-G(y))^2\,dx\,dy
 \ge\frac{1}{12}\left(\frac{a^2}{4C_3}\right)^2=c'a^4.
\]
Thus $\gamma(\tau)\ge c\int_JG^2\,d\xi\ge c'a^4.$ 
\end{proof}

\begin{lemma}
\label{lem:max-evolution}
For sufficiently negative $\tau$, 
\[
 \varrho_{\max}'(\tau)\ge\tfrac15\varrho_{\max}(\tau)
                                      -C\gamma(\tau)^{1/4}.
\]
\end{lemma}
\begin{proof}
Let $m(\tau)=\max G(\xi,\tau)=\varrho_{\max}(\tau)-e^{\tau/4}$. At the maximum at $\xi_\tau^*$, \eqref{eq:decomposition-of-errors} gives 
\[
Q_{1a}=-\frac{m^2}{2(\sqrt{2}+m)},\qquad
Q_{1b}=0,\qquad
Q_2=-\frac{\mathcal E^{\rm orb}(\xi_\tau^*,\tau)}{(\sqrt{2}+m)}.
\]
Since $\varrho_{\max}=e^{\tau/4}+m$ \eqref{eq:G-equation} gives 
\[
{
\varrho_{\max}'
=\varrho_{\max}-\frac34e^{\tau/4}
+G_{\xi\xi}(\xi_\tau^*,\tau)
-\frac{m^2}{2(\sqrt{2}+m)}
-\frac{\mathcal E^{\rm orb}(\xi_\tau^*,\tau)}{(\sqrt{2}+m)}.
}
\]
Using Lemma~\ref{lem:scaled-error} and Lemma~\ref{lem:max-second}, we obtain 
\[
 \varrho_{\max}'\ge\varrho_{\max}-\tfrac34 e^{\tau/4}
                   -Cm^2-C\gamma^{1/4}-Ce^\tau.
\]
By Lemma~\ref{lem:rmax}, $m\ge-Ce^\tau$, so $e^{\tau/4}\le\varrho_{\max}+Ce^\tau$ and
$\varrho_{\max}\ge\frac12e^{\tau/4}$. Since $\varrho_{\max}\to0$, we have $m^2+e^\tau=o(\varrho_{\max})$.
Combining these facts gives the assertion. 
\end{proof}

\begin{lemma}
\label{lem:max-step}
For sufficiently negative $\tau$,
$\varrho_{\max}(\tau-1)\le e^{-1/5}\varrho_{\max}(\tau)+CS(\tau)^{1/4}.$
\end{lemma}
\begin{proof}
Integrate Lemma~\ref{lem:max-evolution} on $[\tau-1,\tau]$. 
\end{proof}

\begin{lemma}
\label{lem:augmented-step}
For sufficiently negative $\tau$, 
\[
 \gamma^+(\tau-1)+\varrho_{\max}(\tau-1)^p
 \le e^{-1}\bigl(\gamma^+(\tau)+\varrho_{\max}(\tau)^p\bigr)
                     +C\delta(\tau)^bS(\tau)+T(\tau).
\]
\end{lemma}
\begin{proof}
Since $p/5>1$, Young's inequality in Lemma~\ref{lem:max-step} gives
\[
 \varrho_{\max}(\tau-1)^p
 \le e^{-1}\varrho_{\max}(\tau)^p+CS(\tau)^{p/4}.
\]
Lemma~\ref{lem:higher-derivatives} gives $S(\tau)\le C\delta(\tau)^{1/2}$. Since $(p/4-1)/2>b$, the last term in the above inequality is at most $C\delta(\tau)^bS(\tau)$. Combining this with Lemma~\ref{lem:spectral-step} gives the assertion. 
\end{proof}
Define, as in \cite{MMS23},  
\[
 \Gamma(\tau)=\sup_{\sigma\le\tau}
                 (\gamma(\sigma)+\varrho_{\max}(\sigma)^p),\qquad
 \Gamma^+(\tau)=\sup_{\sigma\le\tau}
                 (\gamma^+(\sigma)+\varrho_{\max}(\sigma)^p),
\]
\[
 \Gamma^0(\tau)=\sup_{\sigma\le\tau}\gamma^0(\sigma),\qquad
 \Gamma^-(\tau)=\sup_{\sigma\le\tau}\gamma^-(\sigma).
\]

\begin{lemma}
\label{lem:envelope-step}
For sufficiently negative $\tau$,
\begin{align*}
 \Gamma^+(\tau-1)&\le e^{-1}\Gamma^+(\tau)+C\delta(\tau)^b\Gamma(\tau),\\
 |\Gamma^0(\tau-1)-\Gamma^0(\tau)|
       &\le C\delta(\tau)^b\Gamma(\tau),\\
 \Gamma^-(\tau-1)&\ge e\Gamma^-(\tau)-C\delta(\tau)^b\Gamma(\tau).
\end{align*}
Moreover $\delta(\tau)^p\le C\Gamma(\tau)$. 
\end{lemma}
\begin{proof}
Recall that $\|G\|_{L^2(-1,1)}^2\le C\gamma(\tau)$ and by Lemma~\ref{lem:higher-derivatives}, $\|G_\xi\|_{L^2(-1,1)}\le C$. Applying Lemma~\ref{lem:local-H1-interpolation} to $G$ on $[-1,1]$ gives $|G(0,\tau)|\le C\gamma(\tau)^{1/4}$ for all sufficiently negative $\tau$. Combining this with $\varrho^p\le\Gamma$, we obtain $\delta\le C\bigl(\Gamma^{1/p}+\Gamma^{1/4}\bigr)
 \le C\Gamma^{1/p},$ since $p>4$ and $\Gamma\to0$. 
Since $\delta\ge ce^{\tau/4}$ and $8-b=p+b$, $T(\tau)\le C\delta(\tau)^{p+b}
 \le C\delta(\tau)^b\Gamma(\tau).$
Since $S(\tau)\le\Gamma(\tau)$, every error term in Lemmas~\ref{lem:spectral-step} and \ref{lem:augmented-step} is thus bounded by $C\delta(\tau)^b\Gamma(\tau)$. For $\sigma\le\tau$, we have $\delta(\sigma)^b\Gamma(\sigma)\le\delta(\tau)^b\Gamma(\tau)$. The claimed estimates now follow from Lemmas~\ref{lem:spectral-step} and \ref{lem:augmented-step}. 
\end{proof}

\begin{corollary}
\label{cor:source}
We have 
\[
 \|\widehat Q(\tau)\|_{\mathcal H}^2
       \le C\delta(\tau)^{1/200}\Gamma(\tau).
\]
\end{corollary}
\begin{proof}
Use \eqref{eq:localized-residual}, $\gamma\le\Gamma$, and $e^{2\tau}\le C\delta^8\le C\delta^{8-p}\Gamma$. 
\end{proof}

\begin{proposition}
\label{prop:dichotomy}
Exactly one of the following alternatives holds:
\[
 \Gamma^++\Gamma^-=o(\Gamma^0),
 \qquad\hbox{or}\qquad
 \Gamma^0+\Gamma^-=o(\Gamma^+),\qquad \tau\to-\infty.
\]

\end{proposition}
\begin{proof}
Combining the argument used in \cite[Proposition 3.23]{ABDS22} with Lemma~\ref{lem:envelope-step}, the proposition follows. 
\end{proof}

\section{Neutral-mode dominance forces curvature decay}
\label{sec:neutral-mode}
Throughout this section assume that as $ \tau\to-\infty$, 
\begin{equation}\label{eq:neutral}
 \Gamma^+(\tau)+\Gamma^-(\tau)=o(\Gamma^0(\tau)).
\end{equation}
We follow the methods of \cite{MMS23} and \cite[Section~5]{ABDS22} closely to show that we have $R(x)\to 0$ when $d(x,o)\to \infty$, providing details to make the exposition self-contained. 

Recall that $\mathcal H=L^2(\mathbb R,d\nu)$, $d\nu=(4\pi)^{-1/2}e^{-\xi^2/4}\,d\xi$, and the notation from Section~\ref{subsec:local-dichotomy}. 
Let $\mathcal D
 :=\{u\in\mathcal H:\partial_\xi u\in\mathcal H\}$ with the norm  
\[
 \|u\|_{\mathcal D}^2
 :=\int_{\mathbb R}(u^2+|\partial_\xi u|^2)\,d\nu,
\]
where $\partial_\xi u$ is understood in the weak sense. Define
\[
 \alpha(\tau)
 :=\langle\widehat G(\cdot,\tau),h_2\rangle_{\mathcal H},
 \qquad
 A(\tau):=\sup_{\sigma\leq\tau}|\alpha(\sigma)|.
\]
Since $P_0\widehat G=\alpha h_2$ and $\|h_2\|_{\mathcal H}=1$, we have $\Gamma^0=A^2$. 
Combining \eqref{eq:neutral} with Lemma~\ref{lem:envelope-step} and Corollary~\ref{cor:source}, we obtain
\begin{equation}\label{eq:neutral-envelopes}
 \Gamma=(1+o(1))A^2,\qquad
 \delta\le CA^{2/p},\qquad
 \|\widehat Q\|_{\mathcal H}
 \le C\delta^{1/400}A.
\end{equation}
Recall that $A\to 0$ as $\tau \to -\infty$. 

\begin{lemma}
\label{lem:neutral-spectral}
As $\tau\to-\infty$, 
\[
 \|P_+\widehat G\|_{\mathcal D}=o(A),\qquad
\|P_-\widehat G\|_{\mathcal D}
 \le C\delta^{1/400}A=o(A).
\]
In particular, $\|\widehat G-\alpha h_2\|_{\mathcal D}=o(A)$. 
\end{lemma}

\begin{proof}
By \eqref{eq:neutral}, $ \|P_+\widehat G(\tau)\|_{\mathcal H}^2
 \le\Gamma^+(\tau)=o(A(\tau)^2).$ 
Since $\mathcal{H}_+$ is finite dimensional $\|P_+\widehat G\|_{\mathcal D}
\le C\|P_+\widehat G\|_{\mathcal H}=o(A).$  
 Put $v=P_-\widehat G$. Since $P_-$ commutes with $\mathcal L$, Duhamel's formula gives, for $\tau_0<\tau$,
\begin{equation}
\label{eq:duhamel-P-G-hat}
 v(\tau)
 =e^{(\tau-\tau_0)\mathcal L}v(\tau_0)
 +\int_{\tau_0}^{\tau}
   e^{(\tau-\sigma)\mathcal L}P_-\widehat Q(\sigma)\,d\sigma.
\end{equation}
The negative eigenvalues of $\mathcal L$ are at most $-1/2$, so $ \|e^{t\mathcal L}P_-\|_{\mathcal H\to\mathcal H}
 \le e^{-t/2}$ for $t\ge 0$. Since $\widehat G(\tau_0)\to0$ in $\mathcal H$, letting $\tau_0\to-\infty$ in \eqref{eq:duhamel-P-G-hat} yields
\[
 v(\tau)
 =\int_0^\infty e^{t\mathcal L}
                   P_-\widehat Q(\tau-t)\,dt.
\]
Using \eqref{eq:neutral-envelopes} and the monotonicity of
$\delta$ and $A$, we obtain
\[
 \|v(\tau)\|_{\mathcal H}
 \le C\delta(\tau)^{1/400}A(\tau)
       \int_0^\infty e^{-t/2}\,dt
 \le C\delta(\tau)^{1/400}A(\tau).
\]
It remains to estimate $ \|v(\tau)\|_{\mathcal D}$. Applying Duhamel's formula gives 
\begin{equation}
\label{eq:duhamel-P-G-hat-1}
 v(\tau)
 =e^{\mathcal L}v(\tau-1)
 +\int_0^1 e^{t\mathcal L}
                  P_-\widehat Q(\tau-t)\,dt.
\end{equation}
For $u\in C_c^\infty(\R)$, integration by parts yields $\|u\|_{\mathcal D}^2=\langle u,(2-\mathcal L)u\rangle_{\mathcal H}$. Writing $u=\sum_k u_kh_k$ and $\mathcal{L}h_k=\lambda_k h_k$, we obtain for $t\in (0,1]$, 
\[
\begin{aligned}
\|e^{t\mathcal L}P_-u\|_{\mathcal D}^2
&=\sum_{\lambda_k<0}(2-\lambda_k)e^{2t\lambda_k}|u_k|^2\le\sup_{\lambda\le-1/2}(2-\lambda)e^{2t\lambda}
\,\|u\|_{\mathcal H}^2.
\end{aligned}
\]
For some $C>0$, $(2-\lambda)t\leq Ce^{-2t\lambda}$ for all $\lambda\leq -1/2$ and $t\in (0,1]$. Combining this, we get $ \|e^{t\mathcal L}P_-\|_{\mathcal H\to\mathcal D}^2\leq Ct^{-1}$ for $t\in (0,1]$. Consequently,
\begin{align*}
 \|v(\tau)\|_{\mathcal D}
 &\le C\|v(\tau-1)\|_{\mathcal H}
   +C\int_0^1 t^{-1/2}
                  \|\widehat Q(\tau-t)\|_{\mathcal H}\,dt\\
 &\le C\delta(\tau)^{1/400}A(\tau)
       \left(1+\int_0^1t^{-1/2}\,dt\right)\le C\delta(\tau)^{1/400}A(\tau).
\end{align*}
The final statement follows because $\widehat G-\alpha h_2=P_+\widehat G+P_-\widehat G$. This proves the lemma. 
\end{proof}

\begin{lemma}
\label{lem:reference-slope}
$|G_\xi(0,\tau)|=o(A(\tau))$ as $\tau\to-\infty$. 
\end{lemma}

\begin{proof}
Recall that if $v=P_-\widehat G$, then, $G_\xi(0,\tau)
=\partial_\xi P_+\widehat G(0,\tau)
+\alpha(\tau)h_2'(0)+v_\xi(0,\tau).$ 
The first term in this equation is $o(A)$ by Lemma~\ref{lem:neutral-spectral}, $h_2'(0)=0$, so it remains to bound the third term. For $u\in\mathcal H$ and $0<t\le1$, set $w=e^{t\mathcal L}P_-u$.
Applying Lemma~\ref{lem:local-H1-interpolation} to $w_\xi$ gives 
\[
|w_\xi(0)|^2
\le C\|w_\xi\|_{\mathcal H}
\bigl(\|w_\xi\|_{\mathcal H}+\|w_{\xi\xi}\|_{\mathcal H}\bigr).
\]
From the proof of Lemma~\ref{lem:neutral-spectral}, we get $\|w_\xi\|_{\mathcal H}\le Ct^{-1/2}\|u\|_{\mathcal H}$. Parseval identity gives
\[
\|w_{\xi\xi}\|_{\mathcal H}^2
=\frac14\sum_{k\ge3}k(k-1)e^{(2-k)t}|u_k|^2
\le Ct^{-2}\|u\|_{\mathcal H}^2,
\qquad 0<t\le1.
\]
Combining these estimates, we get $|\partial_\xi(e^{t\mathcal L}P_-u)(0)|
\le Ct^{-3/4}\|u\|_{\mathcal H}$ for $0<t\le 1$ and $u\in \mathcal{H}$. Equation~\eqref{eq:duhamel-P-G-hat-1} and  \eqref{eq:neutral-envelopes} now gives 
\[
\begin{aligned}
|v_\xi(0,\tau)|
&\le C\|v(\tau-1)\|_{\mathcal H}
+C\int_0^1t^{-3/4}
       \|\widehat Q(\tau-t)\|_{\mathcal H}\,dt\\
&\le C\delta(\tau)^bA(\tau)
       \left(1+\int_0^1t^{-3/4}\,dt\right)\\
&\le C\delta(\tau)^bA(\tau)=o(A(\tau)).
\end{aligned}
\]
Combining these estimates, we obtain the lemma. 
\end{proof}
Recall that $Q=(\partial_\tau-\mathcal L)G$ and $I_\tau=\left[-\delta(\tau)^{-1/100},\,\delta(\tau)^{-1/100}\right]$. 
\begin{lemma}\label{lem:weighted-neutral-bound}
For sufficiently negative $\tau$,
\[
 \int_{I_\tau}(1+|\xi|)^4(G^2+G_\xi^2)\,d\nu
 \le CA(\tau)^2.
\]
\end{lemma}
\begin{proof}
Put $r=\delta(\tau)^{-1/1000}$ and $v_\pm=P_\pm\widehat G$. For sufficiently negative $\tau$, we have $G=\widehat G=\alpha h_2+v_++v_-$ on $[-\delta^{-1/1000},\delta^{-1/1000}]$. 
Using $\dim \mathcal{H}_+<\infty$ and Lemma~\ref{lem:neutral-spectral}, we get 
\[
\begin{aligned}
 \int_{|\xi|\le \delta^{-1/1000}}(1+|\xi|)^4(G^2+G_\xi^2)\,d\nu
 &\le C\left(\alpha^2+\|v_+\|_{\mathcal D}^2
                  +\delta^{-4/1000}\|v_-\|_{\mathcal D}^2\right)=O(A^2).
\end{aligned}
\]
From Lemma~\ref{lem:higher-derivatives}, $| G|+|\partial_\xi G|\le C$ on $I_\tau$. Thus, 
\[
 \int_{\delta^{-1/1000}\leq |\xi|\le \delta^{-1/100}}
 (1+|\xi|)^4(G^2+G_\xi^2)\,d\nu
 \le Ce^{-cr^2}
 =Ce^{-c\delta^{-1/500}}
 =o(A^2),
\]
where the last equality uses $\delta^p\le CA^2$. Combining the estimates proves the lemma.
\end{proof}

\begin{proposition}
\label{prop:neutral-projection}
We have 
\[
 \int_{I_\tau}Qh_2\,d\nu
       =-2\alpha(\tau)^2+o(A(\tau)^2)+O(e^\tau).
\]
\end{proposition}
\begin{proof}
Recall that $Q=Q_{1a}+Q_{1b}+Q_2$ from \eqref{eq:decomposition-of-errors}. Define $\varepsilon(\tau):=\sup_{I_\tau}(|G|+|G_\xi|)\leq C\delta(\tau)^{1/4}=o(1)$. 
%estimate Q1a. 
We have 
$$G_\xi(0,\tau)\int_{I_\tau}G_\xi(\xi,\tau)h_2(\xi)\,d\nu=o(\varepsilon A).$$
Note that $\left|\frac1{\sqrt2+G}-\frac1{\sqrt2}\right|
\le C\varepsilon.$ Combining these estimates with Lemma~\ref{lem:reference-slope},  we obtain 
\[
\begin{aligned}
\left|\int_{I_\tau}
\left(Q_{1a}+\frac{G^2}{2\sqrt2}
+\frac{G_\xi^2}{\sqrt2}\right)h_2\,d\nu\right|
&\le C\varepsilon A^2+C|G_\xi(0)|A=o(A^2).
\end{aligned}
\]
Using $\int_{\R} h_2^3\,d\nu=2\sqrt2,\int_{\R}(h_2')^2h_2\,d\nu=\sqrt2,$ Lemmas~\ref{lem:neutral-spectral} and \ref{lem:weighted-neutral-bound}, we can obtain 
\[
\int_{I_\tau}
\left(-\frac{G^2}{2\sqrt2}
-\frac{G_\xi^2}{\sqrt2}\right)h_2\,d\nu
=-2\alpha^2+o(A^2).
\]
Hence 
$$\int_{I_\tau}Q_{1a}h_2\,d\nu
       =-2\alpha(\tau)^2+o(A^2).$$
%estimate Q1b. 
Define $B(\xi)=\int_0^\xi
 \frac{G_\xi(x)^2}{(\sqrt2+G(x))^2}\,dx.$ Lemma~\ref{lem:primitive-sobolev} shows that
\[
\|B\|_{L^2(I_\tau,d\nu)}
\le C\left(\int_{I_\tau}G_\xi^4\,d\nu\right)^{1/2}
\le C\|G_\xi\|_\infty\|G_\xi\|_{L^2(I_\tau,d\nu)}
\le C\delta^{1/4}A.
\]
Lemma~\ref{lem:weighted-neutral-bound} gives  $\|G_\xi h_2\|_{L^2(I_\tau,d\nu)}\le CA$. Hence
\[
\left|\int_{I_\tau}Q_{1b}h_2\,d\nu\right|
\le2\|G_\xi h_2\|_2\|B\|_2
\le C\delta^{1/4}A^2=o(A^2).
\]
%estimate Q2
Lemma~\ref{lem:scaled-error} gives $|Q_2|\le Ce^\tau(1+|\xi||G_\xi|)$ on $I_\tau$, hence $\int_{I_\tau}Q_2 h_2\,d\nu=O(e^\tau)$. Combining these estimates, we obtain the proposition. 
\end{proof}

\begin{corollary}
\label{cor:alpha-ode}
Almost everywhere for sufficiently negative $\tau$,
$ \alpha'(\tau)=-2\alpha(\tau)^2+o(A(\tau)^2).$
Moreover, $e^\tau=o(A^2)$. 
\end{corollary}
\begin{proof}
Since $\mathcal Lh_2=0$ and $\widehat Q=(\partial_\tau-\mathcal L)\widehat G$, we get $ \alpha'(\tau)
 =\langle\widehat Q(\tau),h_2\rangle_{\mathcal H}$. Hence  \[
\begin{aligned}
\left|\alpha'-\int_{I_\tau}Qh_2\,d\nu\right|
&\le \int_{I_\tau}\Big|(\vartheta-1)Q+\vartheta_\tau G
-2\vartheta_\xi G_\xi-\vartheta_{\xi\xi}G
+\frac{\xi}{2}\vartheta_\xi G
\Big|h_2\,d\nu\\
&\le C\delta^{-1}
\int_{|\xi|\ge \delta^{-1/100}/2}|h_2|\,d\nu\le C\delta^{-1}e^{-c\delta^{-2/100}}.
\end{aligned}
\]
Since $\delta^p\le CA^2$, $\delta^{-1}e^{-c\delta^{-1/50}} \le CA^2\delta^{-(p+1)}e^{-c\delta^{-1/50}}=o(A^2)$. 
Proposition~\ref{prop:neutral-projection} yields $\alpha'=-2\alpha^2+o(A^2)+O(e^\tau).$ It remains to show that $e^\tau=o(A^2)$. 
Fix $\beta\in(0,1)$. Using $\delta(\tau)\to 0$ and Lemma~\ref{lem:envelope-step}, we obtain for sufficiently negative $\tau$, 
\[
 \Gamma^0(\tau-1)
 \ge\Gamma^0(\tau)-C\delta(\tau)^b\Gamma(\tau)
 \ge\bigl(1-C\delta(\tau)^b\bigr)\Gamma^0(\tau)\ge e^{-\beta}\Gamma^0(\tau).
\]
It follows by iteration that $ A(\tau)^2=\Gamma^0(\tau)\ge c_\beta e^{\beta\tau}$. Consequently, $e^\tau 
 \le c_\beta^{-1}A(\tau)^2e^{(1-\beta)\tau}$ showing $e^\tau=o(A^2)$. 
\end{proof}

\begin{corollary}
\label{cor:sublinear-length}
For all sufficiently large $s$, $d_\pm(s)=o(s).$
\end{corollary}
\begin{proof}
Suppose that $A(\tau)>|\alpha(\tau)|$. 
Since $|\alpha(\tau)|\to 0$ as $\tau\to -\infty$, there exists $\tau_0<\tau$ such that $A(\tau)=\sup_{\sigma\leq \tau}|\alpha(\sigma)|=|\alpha(\tau_0)|$. In a neighborhood of $\tau_0$, Corollary~\ref{cor:alpha-ode} gives $\alpha'\le-cA^2<0$, which is a contradiction. Thus $A(\tau)=|\alpha(\tau)|$ for sufficiently negative $\tau$. 
Corollary~\ref{cor:alpha-ode} now shows that $\alpha$ is strictly decreasing. Hence $\alpha(\tau)=-A(\tau)$ for sufficiently negative $\tau$.

By Cauchy--Schwarz and Lemma~\ref{lem:neutral-spectral}, for $b=\pm 1$, 
\[
 \left|\int_0^{b}\partial_\xi(\widehat G-\alpha h_2)\,d\xi\right|
 \le C\|\widehat G-\alpha h_2\|_{\mathcal D}=o(A).
\]
Since $h_2(\pm1)-h_2(0)=1/(2\sqrt2)$, it follows that
\[
 G(0,\tau)-G(\pm1,\tau)
 =\frac{A(\tau)}{2\sqrt2}+o(A(\tau))
 \ge cA(\tau).
\]
Recall that $G(\cdot,\tau):J_\tau\to \R$ where $J_\tau:=(-e^{\tau/2}d_-(e^{-\tau}),e^{\tau/2}d_+(e^{-\tau}))$ and $G=-\sqrt{2}$ at the endpoints. Concavity of $G$ implies 
\[
 G(0,\tau)-G(\pm1,\tau)
 \le\frac{\sqrt2+G(0,\tau)}{e^{\tau/2}d_\pm(e^{-\tau})}.
\]
Consequently, Corollary~\ref{cor:alpha-ode} therefore gives
\[
 \frac{d_\pm(e^{-\tau})}{e^{-\tau}}
 \le \frac{Ce^{\tau/2}}{A(\tau)}
 =C\left(\frac{e^\tau}{A(\tau)^2}\right)^{1/2}
 \to 0.
\]
\end{proof}

\begin{proposition}
\label{prop:tip-decay}
Under \eqref{eq:neutral},
\[
 \lim_{s\to\infty}R(p_s^-)=\lim_{s\to\infty}R(p_s^+)=0.
\]
\end{proposition}
\begin{proof}
If either limit were positive, Proposition~\ref{prop:zero-tip-sublinear} would give $d_\pm(s)\ge cs$ for some $c>0$ and all sufficiently large $s$. This contradicts Corollary~\ref{cor:sublinear-length}. Hence both limits vanish.
\end{proof}

\section{Excluding positive-mode dominance with \texorpdfstring{$\Z_2$}{Z2}-symmetry}

\label{sec:positive-exclusion}
In this section, we additionally assume that $M$ is $\Z_2$-symmetric, as defined in \cite[Definition~2.1]{Lai20}. We first choose a suitable reference point using the $\Z_2$-symmetry. 

Let $\iota$ denote the reflection. By \cite[Lemma~2.6]{MMS23}, $\iota$ preserves each level set of $\rho$. In the equivariant coordinates $(x',x_4)\in\mathbb R^3\times\mathbb R$ for the $\mathbb Z_2\times O(3)$ action, $\iota(x',x_4)=(x',-x_4)$. Thus $\iota$ exchanges the two halves of the axis $\{x'=0\}$ fixed by $SO(3)$. Since $p_s^\pm$ are the intersections of $\Sigma_s$ with these two half-axes, $\iota(p_s^\pm)=p_s^\mp$.

$\iota$ also preserves $\rho$ and commutes with the level-set flow $\chi_s$, see \cite[Lemma~2.6]{MMS23}.
Since $F(\cdot,s)$ is strictly concave, $F(\cdot,s)$ has a unique maximum. 
The action of $\iota$ on $\Sigma_s$ must fix the spherical orbit with the maximum radius. Choose $q\in\Sigma$ such that for some $s'\gg 1$, $\chi_{s'} q$ lies on the spherical orbit where $F(\cdot,s')$ attains a maximum. Since $\iota$ preserves $\rho$, it commutes with $\chi_s$, so $q_s=\chi_s(q)$ remains on the orbit invariant under $\iota$. Lemma~\ref{lem:rmax} gives $r_{\max}(s)^2\sim2s$ and Lemma~\ref{lem:large-orbit-convergence} gives $ R(q_s)r_{\max}(s)^2\to 2,$ and hence $sR(q_s)\to1$. Rescaling along $q_s$ gives the convergence in Proposition~\ref{prop:level-reference}. Thus, we fix this reference point $q$ in this section. 

Our goal is to prove that the dominance of positive mode in Proposition~\ref{prop:dichotomy} is impossible and conclude the following theorem. 
\begin{theorem}
\label{thm:positive-exclusion}
The positive alternative in Proposition~\ref{prop:dichotomy} cannot occur. Consequently, as $\tau\to-\infty$, 
\[
 \Gamma^+(\tau)+\Gamma^-(\tau)=o(\Gamma^0(\tau)).
\]
\end{theorem}
The argument is by contradiction. That is, in this section, we will assume
\begin{equation}\label{eq:positive-dominance-assumption}
 \Gamma^0(\tau)+\Gamma^-(\tau)=o(\Gamma^+(\tau))
 \qquad (\tau\to-\infty),
\end{equation}
and in section~\ref{sec:coupling-contradiction}, we will derive a contradiction. 

\subsection{Slope on the cylindrical region}
From the choice of $q$, 
\begin{equation}\label{eq:even-profile}
 F(-z,s)=F(z,s),\qquad G(-\xi,\tau)=G(\xi,\tau),\qquad
 G_\xi(0,\tau)=0,
\end{equation}
and $d(s):=d_+(s)=d_-(s)$. Recall the even cutoff $\eta$ and the operator $\mathcal L$ from Section~\ref{sec:dichotomy}. 

\begin{proposition}
\label{prop:exponential-gauge}
As $\tau\to-\infty$, $\Gamma(\tau)=O(e^{\tau/2}), \delta(\tau)=O(e^{\tau/16}).$
\end{proposition}
\begin{proof}
Dominance of positive mode implies $\Gamma\le(1+o(1))\Gamma^+$. Lemma~\ref{lem:envelope-step} gives, for sufficiently negative $\tau$,
\[
 \Gamma^+(\tau-1)
 \le\bigl(e^{-1}+C\delta(\tau)^b(1+o(1))\bigr)\Gamma^+(\tau)
 \le e^{-1/2}\Gamma^+(\tau).
\]
An iteration proves $\Gamma^+(\tau)=O(e^{\tau/2})$ and hence the bound for $\Gamma$. The same lemma gives $\delta^p\le C\Gamma$, where $p=8-1/200<8$. Consequently, $\delta(\tau)\le Ce^{\tau/(2p)}\le Ce^{\tau/16}$ for negative $\tau$. 
\end{proof}
Since $G$ and $\eta$ are even, so is $\widehat G$. Lemma~\ref{lem:cutoff} gives $ \partial_\tau\widehat G=\mathcal L\widehat G+\widehat Q$ almost everywhere. Taking square roots in
\eqref{eq:localized-residual} and using Proposition~\ref{prop:exponential-gauge}, we obtain 
\begin{equation}\label{eq:positive-localized-source}
 \|\widehat Q(\tau)\|_{\mathcal H}
 \le Ce^{\omega\tau}\|\widehat G(\tau)\|_{\mathcal H}+Ce^\tau,
 \qquad \omega:=\frac1{3200},
\end{equation}
almost everywhere, since $\delta^{1/200}\le Ce^{\omega\tau}$ and $\exp(-c\delta^{-1/50})=O(e^\tau)$. Note also that $\|\widehat G(\tau)\|_{\mathcal H}
 \le\Gamma(\tau)^{1/2}=O(e^{\tau/4})$. 

\begin{proposition}
\label{prop:positive-sharp-bounds}
For sufficiently negative $\tau$,
\begin{equation}\label{eq:profile-sharp-bounds}
 \|P_+\widehat G(\tau)\|_{\mathcal H}
 \le C(1+|\tau|)e^\tau,\qquad
 \|(P_0+P_-)\widehat G(\tau)\|_{\mathcal H}\le Ce^\tau.
\end{equation}
Moreover, $\|\widehat Q(\tau)\|_{\mathcal H}\le Ce^\tau$ almost everywhere.
\end{proposition}
\begin{proof}
Fix a sufficiently negative $T$. 
Since $P_\bullet$ commutes with $\mathcal{L}$, we have $\partial_\tau (P_\bullet\widehat G)=\mathcal{L}(P_\bullet\widehat G)+P_\bullet\widehat Q$. 
Using $\|\widehat G(\tau)\|_{\mathcal H}=O(e^{\tau/4})$ and $\|\widehat Q(\tau)\|_{\mathcal H}=O(e^{(1/4+\omega)\tau}+e^\tau)$, Duhamel's formula gives 
\begin{align*}
 (P_+\widehat G)(\tau)
 &=e^{\tau-T}(P_+\widehat G)(T)
   -\int_\tau^T e^{\tau-v}(P_+\widehat Q)(v)\,dv,\\
 (P_0\widehat G)(\tau)
 &=\int_{-\infty}^\tau (P_0\widehat Q)(v)\,dv,\\
 (P_-\widehat G)(\tau)
 &=\int_{-\infty}^\tau
       e^{(\tau-v)\mathcal L}(P_-\widehat Q)(v)\,dv.
\end{align*}
Note that $\|e^{u\mathcal L}P_-\|\le e^{-u}$ for $u\ge0$. If $\|\widehat G(\tau)\|_{\mathcal H}=O(e^{a\tau})$ for some $0<a<1$, then
\[
 \|\widehat Q(\tau)\|_{\mathcal H}
 =O(e^{(a+\omega)\tau}+e^\tau).
\]
When $a+\omega<1$, it follows that $ \|(P_0+P_-)\widehat G(\tau)\|_{\mathcal H}\le Ce^{(a+\omega)\tau}$. On the other hand, 
$$ \|P_+\widehat G(\tau)\|_{\mathcal H}\leq C_T e^\tau+Ce^\tau \int_\tau^T e^{(a+\omega-1)v}\,dv\leq C_T e^{(a+\omega)\tau}.$$
Thus, we obtain $\|\widehat G(\tau)\|_{\mathcal H}=O(e^{(a+\omega)\tau})$. 
Starting with $a=1/4$, iterate while $a+\omega<1$. After finitely many steps, we get $\|\widehat G(\tau)\|_{\mathcal H}=O(e^{a\tau})$ for some $a<1$ with $a+\omega\ge1$.
At this stage, $\|\widehat Q(\tau)\|_{\mathcal H}\le Ce^\tau$, so the Duhamel formulas give $\|(P_0+P_-)\widehat G(\tau)\|_{\mathcal H}\le Ce^\tau$ and 
\[
 \|P_+\widehat G(\tau)\|_{\mathcal H}
 \le C_Te^\tau+C(T-\tau)e^\tau
 \le C(1+|\tau|)e^\tau.
\]
This proves \eqref{eq:profile-sharp-bounds}.  
\end{proof}
Set $a_0(\tau):=\langle\widehat G(\tau),h_0\rangle_{\mathcal H}$ so that $P_+\widehat G=a_0h_0=a_0$. 
\begin{lemma}
\label{lem:positive-pointwise-recovery}
For sufficiently negative $\tau$ and
$|\xi|\le\frac12\delta(\tau)^{-1/100}$,
\begin{align}
 |G(\xi,\tau)-a_0(\tau)|
 &\le Ce^\tau e^{\xi^2/8},
 \label{eq:positive-pointwise-value}\\
 |G_\xi(\xi,\tau)|
 &\le Ce^\tau e^{\xi^2/8}.
 \label{eq:positive-pointwise-slope}
\end{align}
\end{lemma}
\begin{proof}
We prove the estimates for $\xi>0$; the case $\xi<0$ is
analogous, and $\xi=0$ follows by continuity. 
For $0<t\leq 1$ and $v\in \mathcal{H}$, set $w=e^{t\mathcal L}v$. The identity $h_k'=\sqrt{k/2}\,h_{k-1}$ and Parseval's formula gives for $j=0,1,2, 0<t\le 1$, 
\[
 \|\partial_\xi^j w\|_{\mathcal H}^2
 \le C\sup_{k\ge0}(1+k)^j e^{(2-k)t}\|v\|_{\mathcal H}^2
 \le Ct^{-j}\|v\|_{\mathcal H}^2.
\]
On $[\xi-1,\xi]$, we have $y^2\le\xi^2+1$, hence $d\nu(y)\ge ce^{-\xi^2/4}\,dy$ on $[\xi-1,\xi]$.  Applying Lemma~\ref{lem:local-H1-interpolation} to $w$ and $w_\xi$ on $[\xi-1,\xi]$ now yields
\begin{equation}
\label{eq:positive-pt-rec-1}
 |w(\xi)|\le Ct^{-1/4}e^{\xi^2/8}\|v\|_{\mathcal H},
 \qquad
 |w_\xi(\xi)|\le Ct^{-3/4}e^{\xi^2/8}\|v\|_{\mathcal H}.
\end{equation}
Put $U=(P_0+P_-)\widehat G$. Duhamel's formula gives
\[
 U(\tau)=e^{\mathcal L}U(\tau-1)
 +\int_0^1e^{t\mathcal L}
       (P_0+P_-)\widehat Q(\tau-t)\,dt.
\]
Combining this with \eqref{eq:positive-pt-rec-1}, we obtain 
\[
\begin{aligned}
 |U(\xi,\tau)|+|U_\xi(\xi,\tau)|
& \le Ce^{\xi^2/8}\|U(\tau-1)\|_{\mathcal H}+Ce^{\xi^2/8}\int_0^1(t^{-1/4}+t^{-3/4})
   \|(P_0+P_-)\widehat Q(\tau-t)\|_{\mathcal H}\,dt.
\end{aligned}
\]
%Note that the differentiation under the integral sign is valid as $t^{-3/4}$ bound, which is integrable on $(0,1)$. 
By Proposition~\ref{prop:positive-sharp-bounds}, $\|U(\tau-1)\|_{\mathcal H}+\|\widehat Q(\tau-t)\|_{\mathcal H}\leq Ce^\tau$ for almost every $t\in (0,1]$. Thus, 
\[
 |U(\xi,\tau)|+|U_\xi(\xi,\tau)|
 \le Ce^\tau e^{\xi^2/8}.
\]
On $|\xi|\le\frac12\delta(\tau)^{-1/100}$, the cutoff equals one. Since $P_+\widehat G=a_0(\tau)$, we have $U=G-a_0$ and $U_\xi=G_\xi$ there. 
\end{proof}

\begin{proposition}
\label{prop:central-radius-expanding-windows}
As $s\to\infty$,
\begin{equation}\label{eq:positive-central-radius}
 F(0,s)^2=2s+O(\log s).
\end{equation}
\end{proposition}
\begin{proof}
Proposition~\ref{prop:positive-sharp-bounds} and Lemma~\ref{lem:positive-pointwise-recovery} at $\xi=0$ give
$|G(0,\tau)|\le C(1+|\tau|)e^\tau$.
Since $F(0,s)=\sqrt{s}(\sqrt2+G(0,-\log s))$,
\[
 F(0,s)^2
 =2s+2\sqrt2\,sG(0,-\log s)+sG(0,-\log s)^2
 =2s+O(\log s),
\]
where the last term is $O((\log s)^2/s)$. 
\end{proof}

\begin{proposition}
\label{prop:positive-log-window}
For every $\epsilon_2\in(0,4)$, the interval $|z|\le \sqrt{(4-\epsilon_2)s\log s}$ lies in the domain of $F(\cdot,s)$ for all sufficiently large $s$, and 
\begin{equation}\label{eq:positive-log-window}
 \sup_{|z|\le
          \sqrt{(4-\epsilon_2)s\log s}}|F_z(z,s)|=o(s^{-1/2}).
\end{equation}
\end{proposition}
\begin{proof}
Recall that $\log s=-\tau$. The stated interval corresponds to
$|\xi|\le\sqrt{(4-\epsilon_2)(-\tau)}$ after rescaling. Proposition~\ref{prop:exponential-gauge} gives $ \frac12\delta(\tau)^{-1/100}\ge ce^{-\tau/1600},$ so this interval lies in the region where the cutoff has value $1$, and hence in the domain of $G$, for sufficiently large $s$. Using $F_z(z,s)=G_\xi(z/\sqrt{s},-\log s)$ and \eqref{eq:positive-pointwise-slope}, we obtain
\[
 \sup_{|z|\le \sqrt{(4-\epsilon_2)s\log s}}|F_z(z,s)|
 \le C\left(1+\sqrt{|\tau|}\right)
       e^{\tau/2+\epsilon_2\tau/8}
 =o(e^{\tau/2})=o(s^{-1/2}).
\]

\end{proof}

\subsection{Slope on the intermediate region}
\label{sec:intermediate-region-slope}
In this subsection, we will derive bounds on the slope function $F_z$. The methods are inspired by \cite[Section~4]{ABDS22}. We begin with the following observation. 
\begin{lemma}
\label{lem:large-orbit-convergence-2}
For every fixed $b>0$, we have $\sup_{\{F\ge b\sqrt s\}}|F_z|\to 0$. 
\end{lemma}
\begin{proof}
Otherwise, there exist $\varepsilon>0$, $s_j\to\infty$, and $z_j$ such that $F(z_j,s_j)\ge b\sqrt{s_j}, |F_z(z_j,s_j)|\ge\varepsilon.$ 
Let $x_j\in \Sigma_{s_j}$ be any point at signed distance $z_j$ from $\chi_{s_j}q$.  
By Lemma~\ref{lem:large-orbit-convergence}, $R(x_j)\leq Cb^{-2}s_j^{-1}\to 0$ and $(\Sigma_{s_j},R(x_j)g,x_j)$ converges to the round cylinder $\sph^2\times \R$ with basepoint at $\sph^2\times \{0\}$ and the $\R$-factor denoted by coordinate $y$. Then, $R(x_j)F^2\to 2$ and $\sqrt{R(x_j)}(z-z_j)\to y$ smoothly on compact sets. 
Since 
$$F_z=\langle \bar \nabla^ g F,\bar \nabla^g z\rangle_g=\langle \bar \nabla^{R(x_j)g} (R(x_j)^{1/2}F),\bar \nabla^{R(x_j)g} (R(x_j)^{1/2}z)\rangle_{R(x_j)g},$$
it follows that $F_z(x_j)\to 0$, a contradiction.
\end{proof}
The next lemma shows, using \cite{MMS23}, that the slope function $F_z$ satisfies a nonlinear differential inequality with respect to variables $s$ and $r=F(z,s)$. 
\begin{lemma}
\label{lem:intrinsic-slope}
On $0<z<d(s)$, where $F_z<0$, let $r=F(z,s)$ and set $p(r,s)=-F_z(z,s)$. For every sufficiently large $s$ and every $0<r<F(0,s)$, 
\begin{equation}\label{eq:intrinsic-subsolution}
 -p_s\le\mathcal N[p]+C_Er^{-4},
\end{equation}
where 
\begin{equation}
\label{eq:N-definition}
 \mathcal N[p]
 :=p^2p_{rr}
 +(1-p^2)\left(\frac{p_r}{r}+\frac p{r^2}\right),
\end{equation}
and $C_E$ is independent of $r$ and $s$. 
\end{lemma}
\begin{proof}
Set $u(r,s)=p(r,s)^2=F_z(z,s)^2$. The computation in \cite[equation (3.21)]{MMS23}, using Lemmas~\ref{lem:error-sign} and~\ref{lem:error-size}, gives
\[
 -u_s\le uu_{rr}-\frac12u_r^2
 +r^{-2}(1-u)(ru_r+2u)
 +C_0\sqrt u\,r^{-4}.
\]
Substituting $u=p^2$ and dividing by $2p>0$ yields \eqref{eq:intrinsic-subsolution}, with $C_E=C_0/2$.
\end{proof}
Recall that $r_*\in (0,1)$ is the constant such that $\psi_a(r_*/a)\ge3/2$ for all sufficiently large $a$. The proof of the following lemma is given in Appendix~\ref{sec:appendix-A}. 
\begin{lemma}
\label{lem:intrinsic-slope-barriers}
Given $C_E>0$, there exist constants $a_0,K,\Lambda,N>r_*$ and a positive, $C^1$, piecewise smooth function $h:[r_*,\infty)\to(0,\infty)$, with $h(y)=y^{-2}$ for $y\ge N$, having the following property. 

For each fixed $a\ge a_0$, put $\ell_a(s)={\sqrt{2s+\Lambda a^2}}/a,$
and define the $C^1$, piecewise smooth function 
\begin{equation}\label{eq:intrinsic-corrected-barrier}
 P_a(r,s)
 :=
 \sqrt{\psi_a\!\left(\frac r{a\ell_a(s)}\right)}
 +\frac K{\ell_a(s)^2}
       h\!\left(\frac r{\ell_a(s)}\right).
\end{equation}
Then, for $s\in [0,\infty)$, and $r\in [r_*\ell_a(s),\left(1+\frac{a^{-2}}{100}\right)a\ell_a(s)]$, we have 
\begin{equation}\label{eq:intrinsic-strict-super}
 -(P_a)_s>\mathcal N[P_a]+C_Er^{-4}.
\end{equation}
Moreover, 
\begin{equation}\label{eq:intrinsic-barrier-boundary}
 P_a(r_*\ell_a(s),s)>1,
\end{equation}
and, whenever $r\in [N\ell_a(s),\left(1+\frac{a^{-2}}{100}\right)a\ell_a(s)]$, 
\begin{equation}\label{eq:intrinsic-barrier-upper}
 P_a(r,s)\le
 \frac{\sqrt2\,\ell_a(s)}r+\frac K{r^2}.
\end{equation}

\end{lemma}
\begin{remark}
The inequality \eqref{eq:intrinsic-strict-super} holds on each subinterval where $P_a$ is smooth and for both one-sided second derivatives at the joining points of those intervals. 
\end{remark}

The following proposition is the analogue of \cite[Proposition 4.2]{ABDS22} (with $\alpha=\frac12$) in our setting, where we use the modified barrier from Lemma~\ref{lem:intrinsic-slope-barriers}. See also \cite[Proposition~2.13]{Bre20}. 
\begin{proposition}
\label{prop:positive-neck-slope}
There exist $C>0$ and $s_0$ such that, for every $s\ge s_0$ and all $-d(s)<z<d(s)$
\begin{equation}\label{eq:uniform-squared-radius-slope}
 F(z,s)|F_z(z,s)|\le C.
\end{equation}
In particular, the unique point $\bar z(s)\in(0,d(s))$ with $F(\bar z(s),s)=\sqrt s$ satisfies
\begin{equation}\label{eq:positive-neck-slope}
 |F_z(\bar z(s),s)|\le\frac C{\sqrt s}.
\end{equation}
Moreover, for $0\le z\le d(s)$, 
\begin{equation}\label{eq:radius-near-linear-boundary}
 0\le F(z,s)^2\le2C(d(s)-z).
\end{equation}
\end{proposition}

\begin{proof}
On $0<z<d(s)$, where $F_z<0$, use $r=F(z,s)$ as a coordinate. Let $C_E$ be the constant from Lemma~\ref{lem:intrinsic-slope} and let $a_0,K,\Lambda,N>0$, $h$, $P_a$ be obtained by Lemma~\ref{lem:intrinsic-slope-barriers}. By \eqref{eq:positive-central-radius}, 
\begin{equation}\label{eq:intrinsic-comparison-inputs}
 F(0,s)^2\sim2s,\qquad F(0,s)^2\le2s+C_0\sqrt s.
\end{equation}
In particular, $\bar z(s)$ exists and is unique for large $s$. Increasing $\Lambda$ if necessary, \eqref{eq:intrinsic-comparison-inputs} gives, if $s\ge s_0,\ a\ge a_0$, 
\begin{equation}\label{eq:intrinsic-domain}
 \frac{F(0,s)}{\sqrt{2s+\Lambda a^2}}
 <1+\frac{a^{-2}}{100}.
\end{equation}
Indeed, Young's inequality bounds $C_0\sqrt s$ by $\eta s/a^2+C_\eta a^2$, where $\eta>0$ is chosen sufficiently small. 

Fix $a\ge a_0$, and let $\mathcal I_a$ be the set of all $s\ge s_0$ such that
\begin{equation}\label{eq:intrinsic-comparison-set}
 |F_z(z,s)|<P_a(F(z,s),s)
 \quad\text{for all }z\in (0,d(s)),
 F(z,s)\ge r_*\ell_a(s).
\end{equation}
Using Lemma~\ref{lem:large-orbit-convergence-2}, there exists $S_a>0$ such that for every $s\geq S_a$, 
\[
\sup_{\{F\ge \frac{\sqrt{2}r_*}{a}\sqrt s\}}|F_z|<\min_{x\in[r_*/a,\,1+a^{-2}/100]}
\sqrt{\psi_a(x)}.
\]
For all $z\in (0,d(s))$ with $F(z,s)\ge r_*\ell_a(s)$, we have 
\[
 F(z,s)\ge\frac{\sqrt2\,r_*}{a}\sqrt s,
 \qquad
 \frac{r_*}{a}
 \le\frac{F(z,s)}{\sqrt{2s+\Lambda a^2}}
 \le\frac{F(0,s)}{\sqrt{2s+\Lambda a^2}}
 <1+\frac{a^{-2}}{100}.
\]
It follows that $|F_z|\le P_a(F(z,s),s).$ Thus, $[S_a,\infty)\subset\mathcal I_a$. 
\begin{claim}
For every $a\geq a_0$, we have $\mathcal I_a=[s_0,\infty)$. 
\end{claim}
\begin{proof}
Put $ s_*:=\sup\bigl([s_0,\infty)\setminus\mathcal I_a\bigr).$ Then, $s_*\in[s_0,S_a]$. By continuity, there is a point $z_*$ at level $s_*$ where equality first occurs: $ |F_z(z_*,s_*)|
 =P_a(F(z_*,s_*),s_*).$ If $F(z_*,s_*)=r_*\ell_a(s_*),$ we would have
\[
|F_z(z_*,s_*)|\le1
<P_a(r_*\ell_a(s_*),s_*),
\]
which is a contradiction. If $F(z_*,s_*)=F(0,s_*)$, we would obtain $F_z(z_*,s_*)=0<P_a(F(z_*,s_*),s_*)$, which is false. Hence, $r_*\ell_a(s_*)<r_c:=F(z_*,s_*)<F(0,s_*).$ Write $p(r,s)=-F_z(z,s)$ where $r=F(z,s)$ so that $p=P_a$ at $(r_c,s_*)$. At $(r_c,s_*)$, we have
\[
 p_r=(P_a)_r,\qquad
 p_{rr}\le(P_a)_{rr},\qquad
 (p-P_a)_s\le0.
\]
But \eqref{eq:intrinsic-subsolution} and \eqref{eq:intrinsic-strict-super} imply, at $(r_c,s_*)$,
\[
 0\le (P_a-p)_s
 <\mathcal N[p]-\mathcal N[P_a]
 =p^2(p-P_a)_{rr}
  +\frac{1-p^2}{r}(p-P_a)_r
 \le0,
\]
which is impossible. This proves the claim. 
\end{proof}
Fix $s\ge s_0$. If $F(z,s)\leq N\sqrt{\Lambda}$, then 
$F|F_z|\le F\le N\sqrt\Lambda.$
Suppose now that $F(z,s)>N\sqrt\Lambda$.
Since $\ell_a(s)\to\sqrt\Lambda$ as $a\to\infty$,
for all sufficiently large $a$ we have $ F(z,s)>N\ell_a(s)\ge r_*\ell_a(s).$
By the claim and \eqref{eq:intrinsic-barrier-upper},
\[
 F(z,s)|F_z(z,s)|
 \le F(z,s)P_a(F(z,s),s)
 \le\sqrt2\,\ell_a(s)+\frac K{F(z,s)}.
\]
Letting $a\to\infty$ gives
\[
 F(z,s)|F_z(z,s)|
 \le\sqrt{2\Lambda}+\frac K{F(z,s)}
 \le\sqrt{2\Lambda}+\frac K{N\sqrt\Lambda}.
\]The estimate \eqref{eq:uniform-squared-radius-slope} follows. 
Evaluating at $\bar z$ gives \eqref{eq:positive-neck-slope}. Finally,
\[
 F(z,s)^2
 =-2\int_z^{d(s)}F(\zeta,s)F_z(\zeta,s)\,d\zeta
 \le2C(d(s)-z),
\]
which proves \eqref{eq:radius-near-linear-boundary}. 
\end{proof}
A first application of Proposition~\ref{prop:positive-neck-slope} is that, assuming that the positive modes dominate, the limit of scalar curvature at tips is positive. In particular, scalar curvature doesn't vanish at infinity. 
\begin{lemma}
\label{lem:positive-forces-positive-tips}
Under \eqref{eq:positive-dominance-assumption},
\begin{equation}\label{eq:common-positive-tip-limit}
 L:=\lim_{s\to\infty}R(p_s^-)
   =\lim_{s\to\infty}R(p_s^+)\in(0,1).
\end{equation}
\end{lemma}
\begin{proof}
Lemma~\ref{lem:poles}, $\Z_2$-symmetry, and Lemma~\ref{lem:exterior} show that $L\in[0,1)$. Recall that $\bar z(s)$ is chosen to be such that $F(\bar z(s),s)=\sqrt{s}$. Since $F_{zz}<0$ and $F_z(0,s)=0$, Proposition~\ref{prop:positive-neck-slope} gives $ 0\le-F_z(z,s)\le\frac C{\sqrt s}$ for $0\le z\le \bar z(s)$. 
Consequently,
\[
 F(0,s)-\sqrt s
 =\int_0^{\bar z(s)}(-F_z(z,s))\,dz
 \le\frac C{\sqrt s}\,\bar z(s).
\]
Since $F(0,s)/\sqrt s\to\sqrt2$, it follows that $d(s)\ge \bar z(s)\ge cs$ for some $c>0$ and all sufficiently large $s$. If $L=0$, Proposition~\ref{prop:zero-tip-sublinear} gives $d(s)=o(s)$, which is a contradiction. Therefore $L>0$.
\end{proof}

\subsection{The contradiction}
\label{sec:coupling-contradiction}
We now obtain a contradiction to the positive mode dominance. 
The proof below considers the ``deficit'' $\mathcal U(z,s)=(F(0,s)^2-F(z,s)^2)/z$ and shows that it tends to zero at $z=\sqrt{s\log s}$ but stays uniformly positive at $z=\delta s$, for a suitable fixed $\delta>0$. We average this deficit between these two scales using the measure $dz/z$. 
The contradiction is then obtained by analyzing the growth of this average.

\begin{proof}[Proof of Theorem~\ref{thm:positive-exclusion}]
Assume \eqref{eq:positive-dominance-assumption}. By Lemma~\ref{lem:positive-forces-positive-tips} and Proposition~\ref{prop:zero-tip-sublinear}, there exist positive constants $L,\lambda$ such that $R(p_s^\pm)\to L\in (0,1)$ and $d(s)/s\to \lambda:=\sqrt{L/(1-L)}$ as $s\to \infty$. 

By Propositions~\ref{prop:positive-neck-slope} and \ref{prop:central-radius-expanding-windows}, there is $C>0$ such that, for all sufficiently large $s$,
\begin{equation}
\label{eq:contradiction-facts}
 F(0,s)^2\sim2s,\qquad
 0\le -FF_z(z,s)\le C,\qquad
 F(z,s)^2\le 2C(d(s)-z).
\end{equation}
For $z\in (0,d(s))$ and all $s>0$, define 
\[
 \mathcal U(z,s):=\frac{F(0,s)^2-F(z,s)^2}{z}.
\]
\begin{claim}
\label{claim:contradiction-claim-1}
For all sufficiently large $s$, the following holds: $0\leq \mathcal U(z,s)\leq 2C$ for $z\in (0,d(s))$. There exists $\delta\in (0,\lambda)$ such that for all sufficiently large $s$, $\mathcal U(\delta s,s)\ge c >0$. Moreover, $\mathcal U(\sqrt{s \log s},s)\to 0$ as $s\to \infty$. 

Fix such a $\delta$ and define $z_+(s)=\delta s$ and  $z_-(s)=\sqrt{s \log s}$. 
\end{claim}
\begin{proof}
Fix $\delta\in(0,\lambda)$ sufficiently close to $\lambda$ that $2C(\lambda-\delta)<1$. Then, 
\[
 \limsup_{s\to\infty}\frac{F(\delta s, s)^2}{s}
 \le2C(\lambda-\delta)<1.
\]
Combining this with \eqref{eq:contradiction-facts}, we get for all sufficiently large $s$, $\mathcal U(\delta s,s)\ge c >0$. We also have 
\begin{equation}\label{eq:deficit-uniform-bound}
 0\le\mathcal U(z,s)
 =-\frac2z\int_0^zF(\zeta,s)F_z(\zeta,s)\,d\zeta.
\end{equation}
Using \eqref{eq:contradiction-facts} again, $\mathcal{U}(z,s)\leq 2C$ for $z\in (0,d(s))$. Using Proposition~\ref{prop:positive-log-window}, 
\[
\mathcal U(\sqrt{s \log s},s)
\le 2F(0,s)\sup_{z\in \left[0,\sqrt{s\log s}\right]}|F_z(z,s)|=o(1).
\]
This completes the proof of the claim. 
\end{proof}
All constants below may depend on this fixed $\delta$. Concavity of $F(\cdot,s)$ gives
\begin{equation}
\label{eq:contradiction-concavity}
\frac{F(0,s)-F(z_+(s),s)}{z_+(s)}\le \frac{F(0,s)}{d(s)}.
\end{equation}
Consequently, there is a fixed $\theta>0$ such that $F(z_+(s),s)\ge\theta F(0,s)$ for all sufficiently large $s$. Using \eqref{eq:contradiction-facts}, for sufficiently large $s$ and $0\le z\le z_+(s)$, we have 
\begin{equation}\label{eq:contradiction-radius-slope}
 c\sqrt s\le F(z,s)\le C\sqrt s,\qquad
 0\le -F_z(z,s)\le C/\sqrt{s}.
\end{equation}
Define for all $x\gg 1$, 
\[
 J(x)=\frac1x\int_{z_-(e^x)}^{z_+(e^x)}
                  \mathcal U(z,e^x)\,\frac{dz}{z}.
\]
\begin{claim}
\label{claim:contradiction-claim-2}
There is a fixed $\gamma<1$ such that, for all large $x$, 
\begin{equation}\label{eq:deficit-average-upper}
 0\le J(x)\le\gamma\mathcal U(z_+(e^x),e^x).
\end{equation}
In particular, $0\leq J(x)\le C$ for all $x\gg 1$.
\end{claim}
\begin{proof}
For $0<z\le z_+(s)$, concavity of $F(\cdot,s)$ shows that 
\[
\frac{F(0,s)-F(z,s)}{z}
\le
\frac{F(0,s)-F(z_+(s),s)}{z_+(s)}.
\]
Consequently, 
\[
 \mathcal U(z,s)
 \le\frac{2F(0,s)}{F(0,s)+F(z_+(s),s)}
           \mathcal U(z_+(s),s)
 \le\frac2{1+\theta}\mathcal U(z_+(s),s).
\]
Since 
\[
 \frac1x\log\frac{z_+(e^x)}{z_-(e^x)}
 =\frac12+\frac{\log\delta-\frac12\log x}{x}
 \to \frac12,
\]
this gives the claim. 
\end{proof}
The Leibniz rule gives
\begin{equation}
\label{eq:contradiction-leibniz-rule}
 xJ'(x)
 =\mathcal U(z_+(e^x),e^x)
  -\left(\frac12+\frac1{2x}\right)\mathcal U(z_-(e^x),e^x)
  -J(x)
  +\int_{z_-(e^x)}^{z_+(e^x)}e^x\partial_s \mathcal{U}(z,e^x)\,\frac{dz}{z}.
\end{equation}
The next claim bounds the final term from below. 
\begin{claim}
\label{claim:contradiction-claim-3}
For all sufficiently large $s$, 
\begin{equation}\label{eq:deficit-integrated-error}
 \int_{z_-}^{z_+}s\partial_s \mathcal{U}(z,s)\,\frac{dz}{z}
 \ge-o(1).
\end{equation}
\end{claim}
\begin{proof}
In this proof, the quantities are evaluated at $(z,s)$ where $0<z\leq z_+(s)$. 
Proposition~\ref{prop:F-equation} shows 
\[
\begin{aligned}
\partial_s(F^2)
={}&-2FF_{zz}+2(1+F_z^2+E^{\rm orb})\\
&+4FF_z\int_0^z\frac{F_z(\zeta,s)^2}{F(\zeta,s)^2}\,d\zeta
-2FF_z\int_0^zE^{\rm rad}(\zeta,s)\,d\zeta.
\end{aligned}
\]
At $z=0$, this reduces to
\[
\frac{d}{ds}F(0,s)^2
=-2F(0,s)F_{zz}(0,s)+2+2E^{\rm orb}(0,s).
\]
As a result, 
\[
\begin{aligned}
\partial_s\bigl(F(0,s)^2-F(z,s)^2\bigr)
={}&-2F(0,s)F_{zz}(0,s)+2FF_{zz}\\
&-2F_z^2-4FF_z\int_0^z\frac{F_z(\zeta,s)^2}{F(\zeta,s)^2}\,d\zeta\\
&+2E^{\rm orb}(0,s)-2E^{\rm orb}(z,s)+2FF_z\int_0^z E^{\rm rad}(\zeta,s)\,d\zeta.
\end{aligned}
\]
Discarding the nonnegative terms
$-2F(0,s)F_{zz}(0,s)$, $2E^{\rm orb}(0,s)$, and
$-4FF_z\int_0^z F_z^2/F^2\,d\zeta$, we obtain
\[
\begin{aligned}
\partial_s\bigl(F(0,s)^2-F(z,s)^2\bigr)
\ge{}&2FF_{zz}-2F_z^2-2E^{\rm orb}(z,s)+2FF_z\int_0^z E^{\rm rad}(\zeta,s)\,d\zeta.
\end{aligned}
\]
By Lemma~\ref{lem:error-size} and \eqref{eq:contradiction-radius-slope},
\[
 F_z^2+|E^{\rm orb}|\le\frac Cs,\qquad
 \int_0^z|E^{\rm rad}|\,d\zeta\le\frac{Cz}{s^2}.
\]
For all sufficiently large $s$ and $0<z\le z_+(s)$, we obtain 
\begin{equation}\label{eq:deficit-differential-inequality}
s\partial_s\mathcal U
 \ge-\frac{2sF}{z}(-F_{zz})-\frac Cz-\frac Cs
.
\end{equation}
Using $F\le C\sqrt s$ and \eqref{eq:contradiction-radius-slope},
\[
 \begin{aligned}
 2s\int_{z_-}^{z_+}\frac{F(-F_{zz})}{z^2}\,dz
 &\le\frac{Cs^{3/2}}{z_-^2}
              \int_{z_-}^{z_+}(-F_{zz})\,dz\le-\frac{Cs^{3/2}}{z_-^2}F_z(z_+,s)
 \le\frac C{\log s}.
 \end{aligned}
\]
The remaining errors are bounded by $C/z_-$ and $C(\log s)/s$, respectively. This shows \eqref{eq:deficit-integrated-error}. 
\end{proof}
Combining \eqref{eq:contradiction-leibniz-rule} with Claims \ref{claim:contradiction-claim-1}, \ref{claim:contradiction-claim-2}, \ref{claim:contradiction-claim-3}, we conclude that
\[
 xJ'(x)
 \ge(1-\gamma)\mathcal U(z_+(e^x),e^x)-o(1)
 \ge(1-\gamma)c-o(1).
\]
Hence $J'(x)\ge (1-\gamma)c/(2x)$ for all sufficiently large $x$, and integration gives $J(x)\to\infty$. 
This contradicts \eqref{eq:deficit-average-upper}. Thus the dominance of positive modes \eqref{eq:positive-dominance-assumption} is impossible. Combining this with Proposition~\ref{prop:dichotomy} proves Theorem~\ref{thm:positive-exclusion}. 
\end{proof}

\section{Properties of the soliton}
In this section, we derive several properties of the soliton and prove the main theorems of this paper. Let $(M^4,g,f)$ be a complete gradient steady Ricci soliton satisfying $\Ric_g+\nabla^2 f=0.$ Assume that $\Ric_g>0$, $\sec_g\geq0, \nabla f(o)=0$, and $(M,g)$ is $\kappa$-noncollapsed on all scales, $(M,g)$ is $\mathbb Z_2\times O(3)$-symmetric (as in \cite[Definition~2.1]{Lai20}). We may assume $(M,g)$ is not the Bryant soliton. By Proposition~\ref{prop:dichotomy} and Theorem~\ref{thm:positive-exclusion}, we have as $\tau\to -\infty$, 
\[
 \Gamma^+(\tau)+\Gamma^-(\tau)=o(\Gamma^0(\tau)).
\]
Thus, we may apply the results from Section~\ref{sec:neutral-mode}. We obtain that the scalar curvature decays at infinity. 
\begin{corollary}
\label{cor:global-decay}
We have 
\[
 \lim_{d_g(o,x)\to\infty}R(x)=0.
\]
\end{corollary}
\begin{proof}
Combining Proposition~\ref{prop:tip-decay} with \cite[Theorem~4.7(3), Lemma~2.5]{HS26}, we obtain the corollary. 
\end{proof}
With Corollary~\ref{cor:global-decay} and dominance of neutral mode, we may now use the results of \cite[Section~5]{MMS23} to obtain asymptotics for the radius function. 
\begin{proposition}
\label{prop:neutral-fine-asymptotics}
For every fixed $L>0$, uniformly on $|z|\le L\sqrt{s}$,
\footnote{This corrects a typo in \cite[Theorem~1.9 (i)]{MMS23}}
\[
 F(z,s)
 =\sqrt{2s}
 -\frac{z^2-2s}{4\sqrt{2s}\log s}
 +o(1)\frac{\sqrt{s}}{\log s}.
\]
Moreover, for every fixed $\theta\in(0,\tfrac12)$,
\[
 F(z,s)^2
 =2s-\frac{z^2}{2\log s}+o(1)s,
\]
uniformly on $\{z:F(z,s)\ge\theta\sqrt{2s}\}$ as $s\to\infty$. Consequently, for every $\varepsilon\in(0,2)$,
\begin{equation}\label{eq:neutral-intermediate-profile}
 \frac{F(x\sqrt{s\log s},s)^2}{s}
 =2-\frac{x^2}{2}+o(1),
\end{equation}
uniformly for $|x|\le2-\varepsilon$. 
\end{proposition}

\begin{proof}
Corollary~\ref{cor:global-decay} gives decay of scalar curvature at infinity. Together with \eqref{eq:neutral} and the previously established estimates on the errors and radius function $F$, this supplies the hypotheses used in \cite[Section~5]{MMS23}. The assertions follow from \cite[Proposition~5.7, Propositions~5.12--5.13,
and Corollary~5.15]{MMS23}, after the change of variables $z=\sqrt{s}\,\xi$ and $\tau=-\log s$. 
\end{proof}

\begin{proposition}
\label{prop:diameter-asymptotics}
We have $ d_\pm(s)=2(1+o(1))\sqrt{s\log s}$ and $\operatorname{diam}(\Sigma_s,g_{\Sigma_s})=4(1+o(1))\sqrt{s\log s}$. 
\end{proposition}
\begin{proof}
This follows by \cite[Proposition~5.16]{MMS23}. The second statement follows by $\operatorname{diam}(\Sigma_s,g_{\Sigma_s})
 = d_-(s)+d_+(s)$. 
\end{proof}

\begin{proposition}
\label{prop:neutral-tip-asymptotics}
For all sufficiently large $s$,
\begin{equation}
\label{eq:R-at-tip-asymptotic}
R(p_s^\pm)=(1+o(1))\frac{\log s}{s}.
\end{equation}
Moreover, 
\begin{equation}\label{eq:neutral-level-curvature-bound}
 \sup_{\Sigma_s}R= (1+o(1))\frac{\log s}{s},
\end{equation}
and 
\begin{equation}
\label{eq:sharp-neutral-level-curvature-bound}
 \limsup_{d_g(x,o)\to\infty}
 \frac{R(x)d_g(x,o)}{\log d_g(x,o)}=1.
\end{equation}
\end{proposition}
\begin{proof}
Since $R\to 0$ at infinity and the neutral mode dominates, we obtain  \eqref{eq:R-at-tip-asymptotic} from \cite[Proposition~5.18]{MMS23}. Fix a constant $C_0$ such that $R\leq C_0F^{-2}$ on $M\cap \{F>0\}$. Suppose $F_{\mathrm B}$ denotes the radius of spheres of symmetry on the Bryant soliton $\Bry^3$ as a function of distance from its tip. Fix $r\gg 1$ such that $F_{\mathrm B}(r)>(100C_0)^{1/2}$. 

Put $a_s=\sqrt{s/\log s}$ so that $a_s^2R(p_s^\pm)=1+o(1).$ If $\widetilde g_s=a_s^{-2}g$, then $(\Sigma_s,\tilde g_s,p_s^\pm)$ converges to $(\Bry^3,\tilde{g},\tilde{x})$, with scalar curvature $1$ at the tip $\tilde{x}$, in the smooth Cheeger--Gromov sense. Smooth convergence gives for all sufficiently large $s$, 
\[
 \sup_{B_{\Sigma_s}(p_s^\pm,ra_s)}a_s^2R\le 1+o(1).
\]
In particular, $R|_{\Sigma_s}\leq (1+o(1))a_s^{-2}$ on $(-d_-(s),-d_-(s)+ra_s)\cup(d_+(s)-ra_s,d_+(s))$. Smooth convergence shows as $s\to \infty$,  
\[
 \sup_{\partial B_{\Sigma_s}(p_s^\pm,ra_s)}
 \left|a_s^{-1}F-F_{\mathrm B}(r)\right|
 \to 0.
\]
We obtain $F\ge \frac{1}{2}F_{\mathrm B}(r)a_s$ on $\partial B_{\Sigma_s}(p_s^\pm,ra_s)$ for all large $s$. 
Concavity of $F$ shows $F(z,s)\ge \frac{1}{2}F_{\mathrm B}(r)a_s$ for $z\in [-d_-(s)+ra_s,d_+(s)-ra_s]$. 
We obtain for $z\in [-d_-(s)+ra_s,d_+(s)-ra_s]$,
$$R|_{\Sigma_s}(z,s)\le C_0F(z,s)^{-2}\le 4 C_0 F_{\mathrm B}(r)^{-2}a_s^{-2}\leq \frac{4}{100}a_s^{-2}.$$ Combining these facts, we obtain \eqref{eq:neutral-level-curvature-bound}. Fix $\varepsilon\in(0,1)$. By Corollary~\ref{cor:global-decay}, outside a compact set, $1\geq |\nabla\rho|\ge1-\varepsilon$ and integrating, we obtain  
\begin{equation}
\label{eq:rho-linear}
 \rho(x)\leq d_g(o,x)\le C_\varepsilon+\frac{\rho(x)}{1-\varepsilon}.
\end{equation}
Thus, $\rho(x)/d_g(o,x)\to 1$ as $d_g(x,o)\to \infty$. \eqref{eq:sharp-neutral-level-curvature-bound} follows by combining these facts with \eqref{eq:R-at-tip-asymptotic}.
\end{proof}

By \cite[Theorem~1.6]{HS26}, $R(x)d_g(x,\mathcal G)\to0$ as $d(x,o)\to\infty$. The following proposition sharpens this estimate. 
\begin{proposition}
\label{prop:sharp-edge-curvature}
Assume \eqref{eq:neutral}, and set $\mathcal G=\mathcal G_-\cup\mathcal G_+$. Then, 
\begin{equation}\label{eq:sharp-edge-curvature-limsup}
 \limsup_{d_g(x,\mathcal G)\to\infty}
 \frac{R(x)d_g(x,\mathcal G)^2}
      {\log d_g(x,\mathcal G)}=8.
\end{equation}
\end{proposition}

\begin{proof}
Let $a_s=\sqrt{s/\log s}$. Concavity of $F$ shows that on $\Sigma_s$, $\frac{F(z,s)}{d(s)-|z|}\ge \frac{F(0,s)}{d(s)}$. 
Given $\epsilon>0$, \cite[Corollary~2.16 (iv)]{MMS23} shows that for some $r\gg 1$, and sufficiently large $s$, $RF^2\leq 2+\epsilon$ for all $z\in [-d_-(s)+ra_s,d_+(s)-ra_s]$. Combining these facts with Proposition~\ref{prop:diameter-asymptotics} and $F(0,s)^2\sim 2s$, we obtain for $z\in [-d_-(s)+ra_s,d_+(s)-ra_s]$, 
\begin{equation}
 \begin{aligned}
R(z,s)(d(s)-|z|)^2& \le (2+\epsilon)\frac{(d(s)-|z|)^2}{F(z,s)^2}\le (2+\epsilon)\frac{d(s)^2}{F(0,s)^2}\le(4+2\epsilon+o(1))(\log s).
 \end{aligned}
\end{equation}
On the other hand, if $z\in [-d_-(s),-d_-(s)+ra_s]\cup [d_+(s)-ra_s,d_+(s)]$, then $R\cdot (d(s)-|z|)^2\leq r^2R \cdot a_s^2=o(\log s)$. Taking $\epsilon\to 0$, we get for sufficiently large $s$, 
\begin{equation}
\label{eq:sharp-intrinsic-edge-bound}
\sup_{\Sigma_s}R(x)d_{\Sigma_s}(x,\{p_s^+,p_s^-\})^2 \le(4+o(1))(\log s).
\end{equation}
Fix $\varepsilon\in(0,\tfrac12)$ and write $s=\rho(x)$. We consider points with $d_g(x,\mathcal{G})\ge 2$. If $d_g(x,\mathcal G)\le s^{1/2-\varepsilon}$, then \eqref{eq:neutral-level-curvature-bound} gives 
\[
 \frac{R(x)d_g(x,\mathcal G)^2}
      {\log d_g(x,\mathcal G)}
 \le
 C\frac{s^{-2\varepsilon}\log s}
        {\log d_g(x,\mathcal G)}=o(1).
\]
Otherwise, $\log d_g(x,\mathcal G)\ge(\tfrac12-\varepsilon)\log s$, and \eqref{eq:sharp-intrinsic-edge-bound} combined with $d_g(x,\mathcal{G})\leq \min(d(p_s^+,x),d(p_s^-,x))$ gives
\[
 \frac{R(x)d_g(x,\mathcal G)^2}
      {\log d_g(x,\mathcal G)}
 \le\frac{4+o(1)}{\tfrac12-\varepsilon}.
\]
Combining these estimates, we have $\frac{R(x)d_g(x,\mathcal G)^2}
      {\log d_g(x,\mathcal G)}\le\frac{4+o(1)}{\tfrac12-\varepsilon}$ as $s\to \infty$. 
Letting $\varepsilon\to 0$ proves the upper bound in \eqref{eq:sharp-edge-curvature-limsup}. 

For the reference point $q$, $q_s=\chi_s(q)$, we have 
\[
 R(q_s)=\frac{1+o(1)}s,\qquad
 \min_\pm d_{\Sigma_s}(q_s,p_s^\pm)
 =2(1+o(1))\sqrt{s\log s}.
\]
In particular, $d_g(q_s,\mathcal G)\le2(1+o(1))\sqrt{s\log s}$. On the other hand, let $\sigma$ be the unit speed curve joining $q_s$ to a point in $\mathcal{G}$ with length $d_g(q_s,\mathcal{G})$. We have for all large $s$, and all $t\in [0,d_g(q_s,\mathcal{G})]$, 
$$\rho(\sigma(t))=s+\int_{\sigma|_{[0,t]}} \langle\nabla \rho,\sigma'\rangle\geq s-3\sqrt{s\log s}$$
Project $\sigma$ onto $\Sigma_{s'}$ with $s'=s-4\sqrt{s\log s}$ by $\tilde \sigma(t)=\chi_{s'}\circ \chi_{\rho(\sigma(t))}^{-1}\sigma(t)$. Similar to the proof of Proposition~\ref{prop:zero-tip-sublinear}, it follows that 
$$d_g(q_s,\mathcal{G})=L(\sigma)\geq L(\tilde \sigma)\geq  \min_\pm d_{\Sigma_{s'}}(q_{s'},p_{s'}^\pm)
 =2(1+o(1))\sqrt{s'\log s'}$$
This proves that $d_g(q_s,\mathcal G)\ge 2(1+o(1))\sqrt{s\log s}$. Hence
\[
 \frac{R(q_s)d_g(q_s,\mathcal G)^2}
      {\log d_g(q_s,\mathcal G)}
 =8+o(1),
\]
which, together with the uniform upper bound, proves \eqref{eq:sharp-edge-curvature-limsup}. 
\end{proof}
\begin{remark}
\label{rem:proof-of-scalar-decay-main}
Combining propositions~\ref{prop:sharp-edge-curvature} and \ref{prop:neutral-tip-asymptotics} proves Theorem~\ref{thm:scalar-decay}. 
\end{remark}

\begin{proposition}
\label{prop:sharp-neutral-growth}
As $r\to\infty$,
\begin{equation}\label{eq:sharp-neutral-ball-volume}
 \operatorname{Vol}_g B_g(o,r)
 =\frac{128\pi}{15}(1+o(1))
       r^{5/2}\sqrt{\log r}.
\end{equation}
\end{proposition}

\begin{proof}
Set $ a_\pm(s)=\frac{d_\pm(s)}{\sqrt{s\log s}}$ and 
\[
 U_s(x)=\frac{F(x\sqrt{s\log s},s)^2}{s},
 \quad -a_-(s)<x<a_+(s),
\]
and extend $U_s$ by zero outside this interval. Then, $a_\pm(s)=2+o(1)$ so that $U_s$ is supported within $[-3,3]$ for all large $s$. \eqref{eq:neutral-intermediate-profile} gives
\[
 U_s(x)=2-\frac{x^2}{2}+o(1),
\]
locally uniformly on $(-2,2)$. Lemma~\ref{lem:rmax} gives $0\le U_s\le2+o(1)$. Since $g_{\Sigma_s}=dz^2+F(z,s)^2g_{\mathbb S^2}$, we have
\[
 \operatorname{Vol}_3(\Sigma_s)
 =4\pi\int_{-d_-(s)}^{d_+(s)}F(z,s)^2\,dz.
\]
Rescaling and using dominated convergence gives 
\[
 \begin{aligned}
 \frac{\operatorname{Vol}_3(\Sigma_s)}
      {s^{3/2}\sqrt{\log s}}
 &=4\pi\int_{\mathbb R}U_s(x)\,dx=4\pi\int_{-2}^{2}\left(2-\frac{x^2}{2}\right)dx+o(1)=\frac{64\pi}{3}+o(1).
 \end{aligned}
\]
By Corollary~\ref{cor:global-decay} and \eqref{eq:normalization}, $|\nabla\rho|^{-1}=1+o(1)$ uniformly on $\Sigma_s$. Writing $\Omega_s=\{\rho\le s\}$, the coarea formula gives
\[
 \begin{aligned}
 \operatorname{Vol}_4(\Omega_s)
 &=\operatorname{Vol}_4(\Omega_{s_0})
   +\int_{s_0}^{s}\int_{\Sigma_t}
                  \frac{1}{|\nabla\rho|}\,d\mu_t\,dt\\
 &=\frac{64\pi}{3}(1+o(1))
       \int_{s_0}^{s}t^{3/2}\sqrt{\log t}\,dt\\
 &=\frac{128\pi}{15}(1+o(1))
       s^{5/2}\sqrt{\log s}.
 \end{aligned}
\]
Fix $\varepsilon\in(0,1)$. Using \eqref{eq:rho-linear}, there exists $C_\varepsilon>0$ such that $\Omega_{(1-\varepsilon)r-C_\varepsilon}\subset B_g(o,r)\subset\Omega_r$ for all sufficiently large $r$. Hence 
\[
 \frac{128\pi}{15}(1-\varepsilon)^{5/2}+o(1)
 \le
 \frac{\operatorname{Vol}_g B_g(o,r)}
      {r^{5/2}\sqrt{\log r}}
 \le \frac{128\pi}{15}+o(1).
\]
Letting $r\to\infty$ and then $\varepsilon\to 0$ proves \eqref{eq:sharp-neutral-ball-volume}.
\end{proof}
\begin{remark}
\label{rem:proof-of-vol-diam-growth-main}
Corollary~\ref{cor:global-decay} and \cite[Theorem~1.6]{HS26} prove that the asymptotic cone of $(M,g)$ is a ray. Propositions~\ref{prop:neutral-fine-asymptotics}, \ref{prop:diameter-asymptotics}, \ref{prop:sharp-neutral-growth} complete the proof of Theorem~\ref{thm:asymptotic-geometry}.
\end{remark}

We now proceed towards proving Theorem~\ref{thm:dynamical-stability}.
\begin{lemma}
\label{lem:positive-curvature-operator}
We have $\Rm_g>0$ on $M$. 
\end{lemma}
\begin{proof}
Let $p\in M\setminus \mathcal{G}$. By $O(3)$ symmetry, there exists a sphere of symmetry passing through $p$. Choose an orthonormal basis $\{e_1,e_2,e_3,e_4\}$ where $e_1,e_2$ tangent to the sphere and $e_3,e_4$ normal to it. The group of isometries fixing $p$ contains a copy of $O(2)$ acting on $\operatorname{span}\{e_1,e_2\}$ and fixing $e_3,e_4$. This forces the curvature operator $\mathcal{R}:\Lambda^2 T_p M\to \Lambda^2 T_p M$ in the basis
\[
e_1\wedge e_2,\ e_3\wedge e_4,\
e_1\wedge e_3,\ e_1\wedge e_4,\
e_2\wedge e_3,\ e_2\wedge e_4,
\]
to have the block form $\mathcal R=\operatorname{diag}(\kappa_{\text{orb}},\kappa_{\text{nor}},A,B),$
where $A,B$ are symmetric $2\times 2$ matrices. For the same reason, we also have, for $i,j\in\{3,4\}$, 
$\Rm(e_1\wedge e_i,e_1\wedge e_j)
=
\Rm(e_2\wedge e_i,e_2\wedge e_j),$
hence $A=B$. $\kappa_{\rm orb},\kappa_{\rm nor}\ge0$ as $\sec\geq 0$ on $M$. Moreover, for every $v=v_3e_3+v_4e_4\neq 0$, 
\[
(v_3,v_4)A(v_3,v_4)^T
=\big\langle\mathcal R(e_1\wedge v),e_1\wedge v\big\rangle
=|v|^2\sec(e_1,v)\ge0.
\]
Thus $A\ge0$, so $\mathcal R\ge0$ at $p$. Thus, $\Rm_g\geq 0$ on $M\setminus \mathcal{G}$. By continuity $\Rm_g\geq 0$ on $M$. By Hamilton's strong maximum principle, we obtain $\Rm>0$ on $M$. 
\end{proof}

\begin{lemma}
\label{lem:Ric-lower-bound}
For any $a\in(0,1)$, we have $ e^{as}\min_{\Sigma_s}\lambda_{\min}(\Ric_g)\to \infty$ as $s\to \infty$. 
\end{lemma}

\begin{proof}
Write $\lambda=\lambda_{\min}(\Ric_g)$. The soliton identities give
\[
(\Delta \Ric+\langle\nabla\rho,\nabla  \Ric\rangle)_{ij}
=-2\sum_{k,l=1}^4 R_{iklj}R_{kl}.
\]
At any point, choose an orthonormal Ricci eigenbasis
$\{e_i\}_{i=1}^4$ with eigenvalues $\lambda=\lambda_1\le\lambda_2\le\lambda_3\le\lambda_4$. The argument of \cite[Section~2, equations~(3)--(5)]{MW17} yields
\[
\Delta\lambda+\langle \nabla \rho,\nabla \lambda\rangle
\le -2\sum_{i=2}^4 R_{1ii1}\lambda_i
=-2\sum_{i=2}^4\sec_g(e_1, e_i)\lambda_i
\le0,
\]
in the barrier sense, since $\sec_g\ge0$ and $\Ric_g>0$. There exists $L\ge1$ such that $|\nabla\Ric|\le L$ on $M$. 

Fix $y\in M$ and set $r_0=\lambda(y)/(2L)$. 
Then $\lambda\ge\lambda(y)/2$ on $B_g(y,r_0)$. 
Define for $t\geq r_0$, 
\[
\psi(t)=\left(\frac{r_0}{t}\right)^2e^{-2(t-r_0)},
\]
and set $w(x):=\psi(d_g(y,x))$ for $d_g(y,x)\geq r_0$. Then, a direct computation yields $\psi'<0$ and $\psi''+(3t^{-1}+1)\psi'=2\psi$. Write $d=d_g(y,\cdot)$. Since $\Ric\ge0$, we have $\Delta d\le \frac{3}{d}$, $|\nabla d|=1$, and $\langle \nabla \rho,\nabla d\rangle\leq 1$ outside $\operatorname{Cut}(y)\cup\{y\}$. Thus, on $M\setminus (B_g(y,r_0)\cup \operatorname{Cut}(y))$, 
\[
\begin{aligned}
\Delta w+\langle \nabla \rho,\nabla w\rangle
&=\psi''(d)+\psi'(d)\Delta d
  +\psi'(d)\langle\nabla\rho,\nabla d\rangle\\
&\ge\psi''(d)+\left(\frac3d+1\right)\psi'(d)
=2\psi(d)>0.
\end{aligned}
\]
One may use Calabi's trick to have this inequality in the barrier sense on $M\setminus B_g(y,r_0)$. On $\partial B_g(y,r_0)$ we have $w= 1$. 

Consider the function $v(x):=\lambda(x)-\frac{1}{2}\lambda (y)w(x)$ outside $B_g(y,r_0)$. Then $\Delta v+\langle \nabla \rho,\nabla v\rangle<0$ on $M\setminus B_g(y,r_0)$ while $v\geq 0$ on $\partial B_g(y,r_0)$ and $\liminf_{d_g(y,x)\to \infty}v(x)\geq 0$. By the maximum principle it follows that $v\geq 0$ on $M\setminus B_g(y,r_0)$. Thus, we obtain 
\[
\lambda(x)\ge \frac{1}{8L^2}\,\lambda(y)^3e^{-4d_g(x,y)}
\qquad\text{if }d_g(x,y)\ge r_0.
\]
The same estimate holds inside $B_g(y,r_0)$, since
$\lambda(x)\ge\lambda(y)/2$ and $\lambda(y)\le R(y)\le1$.
Hence, for all $x,y\in M$,
\begin{equation}\label{eq:ricci-eigenvalue-comparison}
\lambda(x)\ge \frac{1}{8L^2}\,\lambda(y)^3e^{-4d_g(x,y)}.
\end{equation}
Fix some sufficiently large $s$. If $\lambda(p_s^+)=\Ric(V_s,V_s)$ for some unit vector $V_s\in T_{p_s^+}M$, then write $V_s=\langle \nu,V_s\rangle \nu+W_s$ where $W_s\in T_{p_s^+}\Sigma_s$ and $\nu=\nabla \rho/|\nabla \rho|$. By rotational symmetry, $\lambda(p_s^+)=\langle V_s,\nu\rangle^2 \Ric(\nu,\nu)+\Ric(W_s,W_s)$. 
Dimension reduction gives
$\Ric|_{T\Sigma_s}\ge cR\,g|_{T\Sigma_s}$ at $p_s^+$. We have 
\[
\lambda(p_s^+)\le\Ric(\nu,\nu)
\le\frac{|\nabla R|}{2|\nabla\rho|}
\le CR^{3/2}=o(R).
\]
Consequently, $cR|W_s|^2\le\lambda(p_s^+)=o(R),$ so $|W_s|^2\to0$. We conclude that $\lambda(p_s^+)\geq \frac{1}{2}\Ric(\nu,\nu)$. 

Fix $\varepsilon\in(0,1),\varepsilon\ll a$. Proposition~\ref{prop:neutral-tip-asymptotics} gives, when $s$ is sufficiently large,
\[
R(p_s^+)-R(p_{(1+\varepsilon)s}^+)
=\left(\frac{\varepsilon}{1+\varepsilon}+o(1)\right)
  \frac{\log s}{s}.
\]
Together with
$\frac{d}{ds}R(p_s^+)=-2\Ric(\nu,\nu)(p_s^+)$,
the mean value theorem gives some
$\sigma\in[s,(1+\varepsilon)s]$ such that $\Ric(\nu,\nu)(p_\sigma^+)\ge c_\varepsilon s^{-2}.$ 
Since $\operatorname{diam}_g(\Sigma_s)=o(s)$, we get for every $x\in\Sigma_s$,
\[
d_g(x,p_\sigma^+)
\le \operatorname{diam}_g(\Sigma_s)+d_g(p_s^+,p_\sigma^+)
\le \operatorname{diam}_g(\Sigma_s)+C_0\varepsilon s\leq C'\varepsilon s.
\]
We obtain from \eqref{eq:ricci-eigenvalue-comparison}, 
$$e^{as}\min_{\Sigma_s}\lambda \geq C_L e^{(a-C'\varepsilon)s}s^{-6}$$
Thus, $e^{as}\min_{\Sigma_s}\lambda \to +\infty$ as $s\to \infty$. 
\end{proof}

\begin{proposition}
\label{prop:linear-stability}
The weighted Lichnerowicz operator
\[
 L=-\Delta-\nabla_{\nabla\rho}-2\Rm*,
 \qquad
 d\mu_\rho=e^\rho\,d\mu_g,
\]
has strictly positive bottom of the spectrum on $L^2(d\mu_\rho;S^2T^*M)$. 
\end{proposition}
\begin{proof}
Using Corollary~\ref{cor:global-decay}, Lemmas~\ref{lem:positive-curvature-operator} and~\ref{lem:Ric-lower-bound}, we can apply \cite[Theorem~1.5]{Deruelle} to obtain linear stability. 
\end{proof}
We end this section with the proof of dynamic stability. 
\begin{proof}[Proof of Theorem~\ref{thm:dynamical-stability}]
\cite[Corollary~1.7]{Deruelle} proves dynamic stability of the Bryant soliton. 
Using Corollary~\ref{cor:global-decay}, Lemmas~\ref{lem:positive-curvature-operator} and~\ref{lem:Ric-lower-bound},  \cite[Theorems~1.2--1.3]{Deruelle} yield Theorem~\ref{thm:dynamical-stability}. 
\end{proof}

\appendix

\section{Construction of the barrier function}
\label{sec:appendix-A}

\setcounter{theorem}{0}
\renewcommand{\thetheorem}{\thesection.\arabic{theorem}}
We construct the barrier $P_a$ in Lemma~\ref{lem:intrinsic-slope-barriers} which is used in proving Proposition~\ref{prop:positive-neck-slope}. 
\begin{proof}[Proof of Lemma~\ref{lem:intrinsic-slope-barriers}]
Recall the definition of $\mathcal{N}$ from \eqref{eq:N-definition}. 
Recall from \cite[Theorem~2.1]{Bre20} the singular Bryant profile function $\varphi:(0,\infty)\to(0,\infty)$. Writing $P_0:=\sqrt{\varphi}$ in the ODE of \cite[Theorem~2.1]{Bre20} gives $2P_0\mathcal N[P_0]=0$, hence $\mathcal N[P_0]=0$. 

The barriers $\psi_a$ are constructed using $\varphi$, see \cite[Definition~2.6]{Bre20}. The properties of $\psi_a$ are recorded in \cite[Proposition~2.19]{MMS23}. For some $r_*\in (0,1),N_0\gg 1,a\gg 1$, $\psi_a$ are $C^1$ on $[r_*a^{-1},1+a^{-2}/100]$ and smooth on $(r_*a^{-1},N_0a^{-1})\cup (N_0a^{-1},1+a^{-2}/100)$. 

We define $V_a(\cdot):=\sqrt{\psi_a(\cdot/a)}$ on $[r_*,a+a^{-1}/100]$. Let $\Lambda$ be a positive constant to be chosen and recall that $\ell_a(s):=\sqrt{2s+\Lambda a^2}/a$. 

The parabolic version of $\psi_a$ is defined by $\Psi_a(r,t)=\psi_a(r/\sqrt{-2t})$, following \cite[Corollary~2.8]{Bre20}. After a change of variables $t=-s-\Lambda a^2/2$, we get  
\[
-\partial_s\Psi_a
>
\Psi_a\,\partial_{rr}\Psi_a
-\frac12(\partial_r\Psi_a)^2
+\frac{1-\Psi_a}{r^2}
  \bigl(r\partial_r\Psi_a+2\Psi_a\bigr),
\]
on $r_*\ell_a(s)\le r\le
\left(1+\frac{a^{-2}}{100}\right)a\ell_a(s).$ Since $\Psi_a\left(r,-s-\frac{\Lambda a^2}{2}\right)
=V_a\left(\frac r{\ell_a(s)}\right)^2$,  a direct computation yields
\begin{equation}
\label{eq:V-a-1}
-\partial_s V_a\!\left(\frac r{\ell_a(s)}\right)
-\mathcal N\!\left[V_a\!\left(\frac r{\ell_a(s)}\right)\right]>0,
\end{equation}
on $(r_*\ell_a(s),N_0\ell_a(s))\cup (N_0\ell_a(s),\left(1+\frac{a^{-2}}{100}\right)a\ell_a(s))$. For $x\in (N_0a^{-1},1+\frac{a^{-2}}{100})$, \cite[Definition~2.6,Theorem~2.1]{Bre20} give
\[
 \psi_a'(x)=-2a^{-2}x^{-3}+O(a^{-4}x^{-5}),
 \qquad |\psi_a''(x)|\le Ca^{-2}x^{-4}.
\]
Since $ax\ge N_0$, the first estimate gives $|\psi_a'(x)|\ge ca^{-2}x^{-3}$ for $N_0$ sufficiently large. Differentiating, we obtain
\[
\begin{aligned}
 \frac{(\partial_x^2-x^{-1}\partial_x-x^{-2})\sqrt{\psi_a(x)}}
      {x^{-2}\sqrt{\psi_a(x)}}
 &=\frac{x^2\psi_a''-x\psi_a'}{2\psi_a}
   -\frac{x^2(\psi_a')^2}{4\psi_a^2}-1\le C\frac{a^{-2}x^{-2}}{\psi_a(x)}
   -c\left(\frac{a^{-2}x^{-2}}{\psi_a(x)}\right)^2
 \le C.
\end{aligned}
\]
Substituting $x=r/(a\ell_a(s))$ and using $V_a(r/\ell_a(s))=\sqrt{\psi_a(x)}$ gives
\begin{equation}\label{eq:V-a-2}
 \frac{V_a''\!\left(\frac r{\ell_a(s)}\right)}{\ell_a(s)^2}
 -\frac{V_a'\!\left(\frac r{\ell_a(s)}\right)}{r\ell_a(s)}
 -\frac{V_a\!\left(\frac r{\ell_a(s)}\right)}{r^2}
 \le \frac{C}{r^2}V_a\!\left(\frac r{\ell_a(s)}\right),
\end{equation}
for $r\in (N_0\ell_a(s),\left(1+\frac{a^{-2}}{100}\right)a\ell_a(s))$. \cite[Definition~2.6]{Bre20}, together with the bounds on $\beta_a$ and $\zeta$ in \cite[Section~2]{Bre20}, give $\psi_a(y/a)=\varphi(y)+O(a^{-1})$ in piecewise $C^2([r_*,N])$ as $a\to \infty$.  Hence, $V_a=P_0+O(a^{-1})$. 

With $C_E$ and the above barriers fixed, we will choose the constants $N,D,k_0,k_1,c_h,K,a_0,\Lambda$ in this order, each depending only on the fixed data and the constants to its left. 

Define the $C^1$ function $h:[r_*,\infty)\to (0,\infty)$ by 
\[
 h(y)=
 \begin{cases}
 -yP_0'(y)\left[k_0-\dfrac{k_1}{D}\bigl(e^{D(y-N)}-1\bigr)\right],
      &r_*\le y\le N,\\[4pt]
 y^{-2},&y\ge N.
 \end{cases}
\]
Set $j(y)=-yP_0'(y)>0$ and choose
\[
 k_0=\frac{N^{-2}}{j(N)},\qquad
 k_1=\frac{2N^{-3}+j'(N)k_0}{j(N)}.
\]
The asymptotics for $\varphi$  give $k_1=N^{-2}+O(N^{-4})>0$ for large $N$. These choices give $h(N)=N^{-2}$ and $h'(N)=-2N^{-3}$, and $h$ is $C^1$ on $(r_*,\infty)$. 

Recall that $P_a(r,s)=V_a(\frac{r}{\ell_a(s)})+\frac{K}{\ell^2_a(s)}h(\frac{r}{\ell_a(s)})$ which is smooth on $I:=(r_*\ell_a(s),N_0\ell_a(s))\cup (N_0\ell_a(s) ,N\ell_a(s))\cup (N\ell_a(s),\left(1+\frac{a^{-2}}{100}\right)a\ell_a(s))$. The goal is to show that $-(P_a)_s-\mathcal N[P_a]>C_Er^{-4}$ on $I$. 

\textbf{Case 1. }Suppose $r\in (N\ell_a(s),\left(1+\frac{a^{-2}}{100}\right)a\ell_a(s))$. Substituting $h(y)=y^{-2}$ gives $P_a(r,s)=V_a\!\left(\frac r{\ell_a(s)}\right)+\frac K{r^2}.$ A direct computation yields
\[
\begin{aligned}
-(P_a)_s-\mathcal N[P_a]
={}&-\partial_s V_a\!\left(\frac r{\ell_a(s)}\right)
     -\mathcal N\!\left[V_a\!\left(\frac r{\ell_a(s)}\right)\right]
     +\frac K{r^4}(1-7P_a^2)\\
&-\left(\frac{2K}{r^2}V_a\!\left(\frac r{\ell_a(s)}\right)
          +\frac{K^2}{r^4}\right)
 \left[
 \frac{V_a''\!\left(\frac r{\ell_a(s)}\right)}{\ell_a(s)^2}
 -\frac{V_a'\!\left(\frac r{\ell_a(s)}\right)}{r\ell_a(s)}
 -\frac{V_a\!\left(\frac r{\ell_a(s)}\right)}{r^2}
 \right].
\end{aligned}
\]
Using \eqref{eq:V-a-1}, \eqref{eq:V-a-2}, 
\[
-(P_a)_s-\mathcal N[P_a]
>
\frac K{r^4}\left[
1-7P_a^2
-C\left(
2V_a\!\left(\frac r{\ell_a(s)}\right)^2
+\frac K{r^2}V_a\!\left(\frac r{\ell_a(s)}\right)
\right)
\right].
\]
By \cite[Proposition~2.19]{MMS23}, $\psi_a(x)\le 2a^{-2}x^{-2}$ for $N_0/a\le x\le 1+a^{-2}/100$. For  $r\in (N\ell_a(s),\left(1+\frac{a^{-2}}{100}\right)a\ell_a(s))$, we obtain $V_a(r/\ell_a(s))\le\sqrt2/N$.
Assume $\Lambda\ge K$. Then $K/r^2\le K/(N^2\Lambda)\le N^{-2}$, so $P_a\le \sqrt2/N+N^{-2}$. 
Taking $K>2C_E$ and $N$ sufficiently large, 
\[
-(P_a)_s-\mathcal N[P_a]>\frac{K}{r^4}\left[1-\frac{C}{N^2}\right]>\frac{K}{2r^4}>C_Er^{-4}.
\]

\textbf{Case 2. }Suppose $r/\ell_a(s)\in (r_*,N_0)\cup (N_0,N)$. Recall that 
\[
\mathcal N[W]=W^2W''+(1-W^2)\left(\frac{W'}y+\frac W{y^2}\right).
\]
Differentiating $\mathcal N[P_0+tw]$ at $t=0$ gives the linearization
\[
\mathcal J_0w
=P_0^2w''
+\frac{1-P_0^2}{y}w'
+\left[
\frac{1-P_0^2}{y^2}
+2P_0\left(P_0''-\frac{P_0'}y-\frac{P_0}{y^2}\right)
\right]w.
\]
Differentiating $\mathcal N[P_0(c\,\cdot)](y)=c^2\mathcal N[P_0](cy)=0$ with respect to $c$ at $c=1$ shows that $\mathcal J_0(yP_0')=0$. A direct computation yields
\[
-\mathcal J_0h
=k_1e^{D(y-N)}\bigl(-yP_0'(y)\bigr)P_0(y)^2
\left(D+2\frac{(yP_0')'}{yP_0'}
+\frac{1-P_0^2}{yP_0^2}\right).
\]
Choose $D$ so that the sum in the bracket is at least $1$ on $[r_*,N]$. 
Then, $-\mathcal{J}_0h>c_h>0$ on that interval. 
Now, 
\[
\mathcal N[P_a](r,s)
=\left.\frac1{\ell_a(s)^2}
  \mathcal N\!\left[V_a+\frac K{\ell_a(s)^2}h\right](y)
  \right|_{y=r/\ell_a(s)},
\]
Using $V_a=P_0+O(a^{-1})$ in piecewise $C^2([r_*,N])$ gives 
\[
\mathcal N[P_a]-\mathcal N\!\left[V_a\!\left(\frac r{\ell_a(s)}\right)\right]
=\frac K{\ell_a(s)^4}
\left[
  (\mathcal J_0h)\!\left(\frac r{\ell_a(s)}\right)
  +O\!\left(a^{-1}+\frac K{\ell_a(s)^2}\right)
\right].
\]
Therefore, for $r\in (r_*\ell_a(s),N_0\ell_a(s))\cup (N_0\ell_a(s) ,N\ell_a(s))$, 
\[
\begin{aligned}
-(P_a)_s-\mathcal N[P_a]
={}&-\partial_s V_a\!\left(\frac r{\ell_a(s)}\right)
-\mathcal N\!\left[V_a\!\left(\frac r{\ell_a(s)}\right)\right]+\frac K{\ell_a(s)^4}
\left[-(\mathcal J_0h)\!\left(\frac r{\ell_a(s)}\right)
+O\!\left(a^{-1}+\frac K{\ell_a(s)^2}\right)\right]\\
&>\frac K{\ell_a(s)^4}
\left[c_h+O\!\left(a^{-1}+\frac K{\Lambda}\right)\right],
\end{aligned}
\]
where we also used \eqref{eq:V-a-1},  $-\mathcal{J}_0h>c_h>0$, and $\ell_a(s)^2\ge\Lambda$. Choosing $a_0$ and $\Lambda/K$ sufficiently large makes the term in the last bracket at least $c_h/2$. Therefore
\[
-(P_a)_s-\mathcal N[P_a]
>\frac{c_hK}{2\ell_a(s)^4}
\ge\frac{c_hKr_*^4}{2r^4}
>C_Er^{-4},
\]
where the last inequality follows from choosing $K>2C_E/(c_hr_*^4)$. 

Combining Cases 1 and 2 (and the choice of constants determined there), we obtain $-(P_a)_s-\mathcal N[P_a]>C_Er^{-4}$ on $I$. 
$\psi_a(r_*/a)\ge 3/2$ and $h(r_*)>0$ give
\[
 P_a(r_*\ell_a(s),s)
 =\sqrt{\psi_a(r_*/a)}+\frac{K}{\ell_a(s)^2}h(r_*)
 \ge\sqrt{\frac32}>1.
\]
For $r\ge N\ell_a(s)$, we have $h(r/\ell_a(s))=\ell_a(s)^2/r^2$. Since $N\ge N_0$, the bound $\psi_a(x)\le 2a^{-2}x^{-2}$ at $x=r/(a\ell_a(s))$ gives
\[
 P_a(r,s)
 =\sqrt{\psi_a\!\left(\frac{r}{a\ell_a(s)}\right)}+\frac{K}{r^2}
 \le \frac{\sqrt2\,\ell_a(s)}{r}+\frac{K}{r^2}.
\]
This completes the proof of Lemma~\ref{lem:intrinsic-slope-barriers}. 
\end{proof}

\section{Some auxiliary results}
\label{sec:appendix-B}
We will be using the following one-dimensional Gaussian Sobolev inequality in the paper. Recall that $d\nu(\xi)=\frac{1}{\sqrt{4\pi}}e^{-\xi^2/4}\,d\xi.$
\begin{lemma}
\label{lem:primitive-sobolev}
Let $r>0$ and $h\in C^1([-r,r])$ satisfy $h(0)=0$. Then
\[
 \int_{-r}^r h^2\,d\nu
 \le 2\int_{-r}^r |h'|^2\,d\nu.
\]
\end{lemma}
\begin{proof}
Write $d\nu=w(\xi)\,d\xi$, so $w'=-\xi w/2$. Expanding the square and integrating by parts gives
\[
 \int_\epsilon^r |h'|^2\,d\nu
 =\int_\epsilon^r\left(h'-\frac{h}{\xi}\right)^2\,d\nu
  +\frac12\int_\epsilon^r h^2\,d\nu
  +\frac{h(r)^2w(r)}{r}-\frac{h(\epsilon)^2w(\epsilon)}{\epsilon}.
\]
Since $h(0)=0$ and $h(\xi)=O(\xi)$, $h(\epsilon)^2/\epsilon\to 0$ as $\epsilon\to 0$. Taking $\epsilon\to 0$ in the above inequality yields $2\int_0^r |h'|^2\,d\nu\geq \int_0^r h^2\,d\nu$. Adding the corresponding estimate for $h(-\xi)$ gives the result. 
\end{proof}
The following interpolation inequality is often useful in converting $H^1$ norms to $L^\infty$ norms. 
\begin{lemma}
\label{lem:local-H1-interpolation}
Let $I$ be an interval of length $\ell>0$. For every $f\in H^1(I)$,
\[
 \|f\|_{L^\infty(I)}^2
 \le \ell^{-1}\|f\|_{L^2(I)}^2
      +2\|f\|_{L^2(I)}\|f'\|_{L^2(I)}.
\]
\end{lemma}
\begin{proof}
For absolutely continuous $f$, and $x,y\in I$, we have 
\[
|f(x)|^2-|f(y)|^2
=2\int_y^x f(t)f'(t)\,dt
\le 2\int_I |f(t)f'(t)|\,dt.
\]
Average over $y$, apply Cauchy--Schwarz to the right hand side, and take the supremum over $x$ to obtain the inequality.   
\end{proof}


\begin{thebibliography}{ABDS22}

\bibitem[ABDS22]{ABDS22}
S.~Angenent, S.~Brendle, P.~Daskalopoulos, and N.~\v{S}e\v{s}um,
\emph{Unique asymptotics of compact ancient solutions to
three-dimensional Ricci flow},
Comm. Pure Appl. Math. \textbf{75} (2022), no.~5, 1032--1073.
\url{https://doi.org/10.1002/cpa.21955}.

\bibitem[BL26]{BL26}
R.~H.~Bamler and Y.~Lai,
\emph{Cylindrical Ricci flows and the PDE-ODI principle},
arXiv:2611.29246 (2026).
\url{https://arxiv.org/abs/2611.29246}.

\bibitem[Bre20]{Bre20}
S.~Brendle,
\emph{Ancient solutions to the Ricci flow in dimension 3},
Acta Math. \textbf{225} (2020), no.~1, 1--102.
\url{https://doi.org/10.4310/ACTA.2020.v225.n1.a1}.

\bibitem[Bry05]{Bry05}
R.~L.~Bryant,
\emph{Ricci flow solitons in dimension three with $SO(3)$-symmetries},
preprint, Duke University, 2005.
\url{https://sites.math.duke.edu/~bryant/3DRotSymRicciSolitons.pdf}.

\bibitem[CMZ25]{CMZ25}
P.-Y.~Chan, Z.~Ma, and Y.~Zhang,
\emph{Dimension reduction for positively curved steady solitons},
Ann. Global Anal. Geom. \textbf{67} (2025), no.~4, Paper No.~20.
\url{https://doi.org/10.1007/s10455-025-10001-8}.

\bibitem[CZ22]{ChanZhu22}
P.-Y.~Chan and B.~Zhu,
\emph{On a dichotomy of the curvature decay of steady
Ricci solitons},
Adv. Math. \textbf{404} (2022), part B,
Paper No.~108458, 40 pp.
\url{https://doi.org/10.1016/j.aim.2022.108458}.

\bibitem[DZ20]{DZ20}
Y.~Deng and X.~Zhu,
\emph{Classification of gradient steady Ricci solitons with linear curvature decay},
Sci. China Math. \textbf{63} (2020), no.~1, 135--154.
\url{https://doi.org/10.1007/s11425-019-1548-0}.

\bibitem[Der12]{Deruelle12}
A.~Deruelle,
\emph{Steady gradient Ricci soliton with curvature in $L^1$},
Comm. Anal. Geom. \textbf{20} (2012), no.~1, 31--53.
\url{https://doi.org/10.4310/CAG.2012.v20.n1.a2}.

\bibitem[Der15]{Deruelle}
A.~Deruelle,
\emph{Stability of non compact steady and expanding gradient
Ricci solitons},
Calc. Var. Partial Differential Equations
\textbf{54} (2015), no.~2, 2367--2405.
\url{https://doi.org/10.1007/s00526-015-0868-4}.

\bibitem[F26]{Fang2026}
H.~Fang,
\emph{$O(3)$-symmetry of four-dimensional steady Ricci solitons},
arXiv:2610.36563 (2026).
\url{https://arxiv.org/abs/2610.36563}.

\bibitem[Ham95]{Ham95}
R.~S.~Hamilton,
\emph{The formation of singularities in the Ricci flow},
Surv. Differ. Geom. \textbf{2} (1995), 7--136.
\url{https://doi.org/10.4310/SDG.1993.v2.n1.a2}.

\bibitem[HS26]{HS26}
A.~G.~Hebbar and N.~\v{S}e\v{s}um,
\emph{Asymptotic geometry of four-dimensional steady solitons},
arXiv:2603.21397v2 (2026).
\url{https://arxiv.org/abs/2603.21397v2}.

\bibitem[Lai24]{Lai20}
Y.~Lai,
\emph{A family of 3D steady gradient solitons that are flying wings},
J. Differential Geom. \textbf{126} (2024), no.~1, 297--328.
\url{https://doi.org/10.4310/jdg/1707767339}.

\bibitem[MMS26]{MMS23}
Z.~Ma, H.~Mahmoudian, and N.~\v{S}e\v{s}um,
\emph{Unique asymptotics of steady Ricci solitons with symmetry},
J. Eur. Math. Soc. (2026), published online first.
\url{https://doi.org/10.4171/JEMS/1774}.

\bibitem[MSW19]{MSW19}
O.~Munteanu, C.-J.~A.~Sung, and J.~Wang,
\emph{Poisson equation on complete manifolds},
Adv. Math. \textbf{348} (2019), 81--145.
\url{https://doi.org/10.1016/j.aim.2019.03.019}.

\bibitem[MW17]{MW17}
O.~Munteanu and J.~Wang,
\emph{Positively curved shrinking Ricci solitons are compact},
J. Differential Geom. \textbf{106} (2017), no.~3, 499--505.
\url{https://doi.org/10.4310/jdg/1500084024}.

\end{thebibliography}
\end{document}